\documentstyle[12pt]{article}

\catcode`\@=11

\ifcase\@ptsize
 \font\tenmsa=msam10
 \font\sevenmsa=msam7
 \font\fivemsa=msam5
 \font\tenmsb=msbm10
 \font\sevenmsb=msbm7
 \font\fivemsb=msbm5
 \font\teneuf=eufm10
 \font\seveneuf=eufm7
 \font\fiveeuf=eufm5
\or
 \font\tenmsa=msam10 scaled \magstephalf
 \font\sevenmsa=msam7 scaled \magstephalf 
 \font\fivemsa=msam5 scaled \magstephalf  
 \font\tenmsb=msbm10 scaled \magstephalf
 \font\sevenmsb=msbm7 scaled \magstephalf 
 \font\fivemsb=msbm5 scaled \magstephalf  
 \font\teneuf=eufm10 scaled \magstephalf
 \font\seveneuf=eufm7 scaled \magstephalf
 \font\fiveeuf=eufm5 scaled \magstephalf
\or
 \font\tenmsa=msam10 scaled \magstep1
 \font\sevenmsa=msam7 scaled \magstep1 
 \font\fivemsa=msam5 scaled \magstep1  
 \font\tenmsb=msbm10 scaled \magstep1
 \font\sevenmsb=msbm7 scaled \magstep1 
 \font\fivemsb=msbm5 scaled \magstep1  
 \font\teneuf=eufm10 scaled \magstep1
 \font\seveneuf=eufm7 scaled \magstep1
 \font\fiveeuf=eufm5 scaled \magstep1
\fi

\newfam\msafam
\newfam\msbfam
\newfam\euffam

\textfont\msafam=\tenmsa  \scriptfont\msafam=\sevenmsa
  \scriptscriptfont\msafam=\fivemsa
\textfont\msbfam=\tenmsb  \scriptfont\msbfam=\sevenmsb
  \scriptscriptfont\msbfam=\fivemsb
\textfont\euffam=\teneuf  \scriptfont\euffam=\seveneuf
  \scriptscriptfont\euffam=\fiveeuf

\def\hexnumber@#1{\ifnum#1<10 \number#1\else
 \ifnum#1=10 A\else\ifnum#1=11 B\else\ifnum#1=12 C\else
 \ifnum#1=13 D\else\ifnum#1=14 E\else\ifnum#1=15
F\fi\fi\fi\fi\fi\fi\fi}

\def\Bbb{\ifmmode\let\next\Bbb@\else
\def\next{\errmessage{Use \string\Bbb\space only in math
mode}}\fi\next}
\def\Bbb@#1{{\Bbb@@{#1}}}
\def\Bbb@@#1{\fam\msbfam#1}

\def\frak{\ifmmode\let\next\frak@\else
\def\next{\errmessage{Use \string\frak\space only in math
mode}}\fi\next}
\def\frak@#1{{\frak@@{#1}}}
\def\frak@@#1{\fam\euffam#1}

\def\goth{\ifmmode\let\next\frak@\else
\def\next{\errmessage{Use \string\goth\space only in math
mode}}\fi\next}

\def\mapdown#1{{\Big\downarrow\rlap{$\vcenter {\hbox {$#1$}}$}}}
\def\mapdownalt#1{{\Big\downarrow\llap{$\vcenter{\hbox {$#1$}}$~~}}}
\def\mapright#1{{\smash{\mathop{\longrightarrow}\limits^{#1}}}}

\def\mapne#1{{\nearrow\llap{$\vcenter {\hbox {$#1$}}$~~~~}}}

\def\sqr#1#2{{\vcenter{\vbox{\hrule height.#2pt
        \hbox{\vrule width.#2pt height#1pt \kern#1pt
           \vrule width.#2pt}
        \hrule height.#2pt}}}}
\def\square{\mathchoice\sqr56\sqr56\sqr{2.1}3\sqr{1.5}3}

\def\N{{\Bbb N}}
\def\Z{{\Bbb Z}}
\def\Q{{\Bbb Q}}
\def\R{{\Bbb R}}
\def\C{{\Bbb C}}
\def\P{{\Bbb P}}

\newtheorem{thm}{Theorem}[section]
\newtheorem{lem}[thm]{Lemma}
\newtheorem{prop}[thm]{Proposition}
\newtheorem{cor}[thm]{Corollary}
\newtheorem{defn}[thm]{Definition}
\newtheorem{pagph}[thm]{($\!$}
\newtheorem{ex}[thm]{Example}
\newtheorem{rk}[thm]{Remark}
\newtheorem{conj}[thm]{Conjecture}

\newenvironment{theorem}{\begin{thm}\sl}{\end{thm}\rm}
\newenvironment{lemma}{\begin{lem}\sl}{\end{lem}\rm}
\newenvironment{proposition}{\begin{prop}\sl}{\end{prop}\rm}
\newenvironment{corollary}{\begin{cor}\sl}{\end{cor}\rm}

\newenvironment{remark}{\begin{rk}\rm}{\end{rk}}

\def\proof{\noindent{\bf Proof.~}}
\def\endproof{\quad $\square$}

\catcode`\@=12

\def\H{{\Bbb H}}
\def\L{{\Bbb L}}
\def\V{{\Bbb V}}
\def\W{{\Bbb W}}

\def\D{{\cal D}}
\def\cF{{\cal F}}
\def\cH{{\cal H}}
\def\cP{{\cal P}}
\def\cV{{\cal V}}
\def\cW{{\cal W}}

\def\gD{{\frak D}}
\def\gP{{\frak P}}
\def\f{{\frak f}}
\def\g{{\frak g}}
\def\p{{\frak p}}

\def\v{\vec{v}}

\def\LQ{{\Lambda\otimes\Q}}

\def\Vbar{\overline{\cV}}
\def\Fbar{\overline{\cF}}
\def\Wbar{\overline{\cW}}
\def\Deltabar{\overline{\Delta}}

\def\w{\omega}
\def\bw{\mbox{\boldmath$\omega$}}
\def\bL{\lambda}                   
\def\O{{\cal O}}
\def\X{\overline{X}}
\def\E{\overline{E}}
\def\nablabar{\overline{\nabla}}
\def\Deltabar{\overline{\Delta}}

\def\Hdel{H_{\cal D}}
\def\Hom{{\rm Hom}}
\def\Hcts{H_{\rm cts}}
\def\End{{\rm End}}
\def\Aut{{\rm Aut}}
\def\Ext{{\rm Ext}}

\def\01{\{0,1\}}
\def\half{1/2}
\def\span{{\rm span}}
\def\spec{{\rm spec\, }}
\def\vol{{\rm vol\, }}
\def\sgn{{\rm sgn}}
\def\im{{\rm im}}
\def\ad{{\rm ad}}

\def\comp{~\widehat{\!}{\;}}
\def\comptensor{\widehat\otimes}

\def\dlog#1{{d#1 \over #1}}
\def\row#1#2#3{#1_{#2} \ldots #1_{#3}}

\begin{document}

\begin{center}
{\LARGE Classical Polylogarithms} \\
\medskip

{\bf March, 1992 \\

Richard M. Hain}\footnote{Supported in part by grants from the
National
Science Foundation.}\\
Department of Mathematics\\
Duke University\\
Durham, NC 27706
\end{center}
\bigskip

\section{Introduction}
\label{intro}

This article is an introduction to classical polylogarithms. After
establishing some of their basic properties, we
present several examples of Spencer Bloch where the dilogarithm is
used to construct the second regulator.  We also construct the
polylogarithm local systems and show that each underlies a Tate
variation of mixed Hodge structure \cite{deligne_letter}. We conclude
by giving an exposition of a motivic description of the polylogarithm
local systems. Most of the results in this paper were discovered by
Bloch, Deligne, Ramakrishnan, Suslin and Beilinson.

Let $k$ be  a positive integer.
The $k$th {\it polylogarithm} $\ln_k x$ is defined by
\begin{equation}\label{series}
\ln_k x = \sum_{n=1}^\infty {x^n \over n^k}.
\end{equation}
This converges in the unit disk to a holomorphic
function.
The first polylogarithm, $\ln_1 x$, is just $- \log(1-x)$. The
second,
$$
\ln_2 x = \sum_{k=1}^\infty  {x^n \over n^2}
$$
is called the {\it dilogarithm}, and was defined by Euler in 1768.
The higher polylogarithms were defined by Spence in 1809
(cf. \cite{lewin}).

It is believed that the $k$th regulator
$$
c_k : K_m(X) \to \Hdel^{2k-m}(X,\Z(k))
$$
from the algebraic $K$-theory of a complex algebraic variety $X$ (and
therefore all varieties of finite type over $\Q$) to its
Deligne cohomology can be expressed in terms of the $k$th
polylogarithm.
If true, this would generalize
the classical fact that the logarithm occurs as the first Chern class
$$
c_1 : K_0(X) \to H^2(X,\Z(1))
$$
and its single valued cousin, $\log|\phantom{x}| :  \C^\ast \to \R$,
occurs as the regulator
$$
c_1 : K_1(\C) \approx \C^\ast \to \R
\approx\Hdel^1({\rm spec}\, \C,\R(1)).
$$

In the case where $X$ is $\spec \C$, this should mean
that some single valued cousin $D_k : \C - \01 \to \R$ of $\ln_k$
should represent a multiple of the Borel regulator element
$$
b_k \in H^{2k-1}(GL_k(\C)^\delta,\R),
$$
the cohomology class which gives rise to the regulator \cite{borel}
(cf. also \cite[\S 2.2]{ramakrishnan_survey}, \cite[\S 7.2]{D-H-Z}.
(Here, $GL_k(\C)^\delta$ denotes the general linear group viewed as
a discrete group.)
The cocycle condition would then be a
functional equation satisfied by $D_k$ which generalizes
the 3-term functional equation satisfied by $D_1 = \log
|\phantom{x}|$.
It is further believed that all of the rational $K$-theory of a field
$F$ should come from $F - \01$, and that the relations should all
correspond to canonical functional equations satisfied by the $D_k$.
Such statements are often referred to as the Zagier conjecture.%
\footnote{Zagier's original conjecture \cite{zagier} asserted that
the
value at the positive integer $m$ of the Dedekind zeta function  of a
number field $F$ could be expressed as a determinant of values of
$D_m$
at $F$ rational points of $\P^1 - \{0,1,\infty\}$.}
For example, for all fields, we may express the familiar fact
$$
K_1(F) = F^\times
$$
as
$$
K_1(F) = \left[\coprod_{x \in F - \01} \Z \right]/{\cal R}
$$
where  the relations $\cal R$ are generated by
$$
[x] - [xy] + [y] = 0\quad \hbox{and} \quad [x] + [x^{-1}] = 0,
$$
where $x,y\in F - \01$ and $xy \neq 1$. These are the analogues of
the functional equations
$$
D_1(x) - D_1(xy) + D_1(y) = 0 \quad \hbox{and}
\quad D_1(x) + D_1(x^{-1}) =  0.
$$
Note also that the functional equation in this case is precisely the
condition that $D_1$ represent an element of $H^1(GL_1(\C),\R)$.
The corresponding story for the dilogarithm has been worked out by
Bloch \cite{bloch-irvine} and Suslin \cite{suslin}.  We give an
account of this story in Section \ref{reg-k3}.

Useful references for basic material in this paper include
\cite{milnor}
and [this volume] for algebraic $K$-theory, \cite{lewin} for a
comprehensive
reference on classical aspects of polylogarithms,
\cite{hain-zucker_2}
for Tate variations of mixed Hodge structure, and \cite{geom} for
basic facts about iterated integrals and the mixed Hodge theory of
the
fundamental group.  The  book \cite{lewin_2} is a useful reference
for more recent developments. (Also, references to other articles in
this book.)

\noindent{\bf Notation:} The group of units of a ring $R$ will be
denoted
by $R^\times$.  When $\Lambda$ is $\Z$, $\Q$, or $\R$, $\Lambda(k)$
will
denote the subgroup $(2\pi i)^k\Lambda$ of $\C$.  It will also be
used to
denote the Hodge structure of type $(-k,-k)$ which has this abelian
group
as its lattice.

\section{Monodromy}
\label{monodromy}

An easy power series manipulation yields the formula
$$
\ln_k x = \int_0^x \ln_{k-1}z\, {dz \over z},
$$
where $x$ lies in the unit disk and $k\ge 2$.
It follows, by an induction argument, that each polylogarithm can be
analytically
continued to a multivalued holomorphic function on $\C-\01$. In this
section we determine the monodromy of the polylogarithms. This was
first
computed by Ramakrishnan in \cite{ramakrishnan_monod}.

Set
$$
\w_0 = {dz \over z} \hbox{ and } \w_1 = {dz \over 1-z}.
$$
Let
$$
\omega = \pmatrix{
0 & \w_1 & 0 & \cdots & 0 \cr
  & \ddots & \w_0 & \ddots & \vdots \cr
\vdots &  & \ddots & \ddots & 0 \cr
  &	&	& \ddots & \w_0 \cr
0 &	& \cdots &	& 0 \cr }
\in H^0(\P^1,\Omega^1_{\P^1}(\log\{0,1,\infty\}))\otimes gl_n(\C).
$$
Consider the first order linear differential equation
\begin{equation}\label{diff-eqn}
d\bL = \bL \bw
\end{equation}
where $\bL$ is a possibly multivalued function $\C - \01 \to
\C^{n+1}$.

Denote the $k$th power of the standard logarithm $\log x = \int_1^x
\w_0$
by $\log^k x$.
Let
$$
\Lambda(x) = \pmatrix{
1 & \ln_1 x & \ln_2 x & \cdots & \cdots & \ln_n x \cr
0 & 2\pi i & 2\pi i\log x & & \cdots & {2\pi i\over n!}\log^{n-1}x
\cr
0 & 0 & (2 \pi i)^2 & \ddots & \cdots & {(2\pi i)^2 \over
(n-1)!}\log^{n-1}x
\cr
\vdots & & \ddots & \ddots &  & \vdots \cr
\vdots & & & 0 & (2\pi i)^{n-1} & (2\pi i)^{n-1}\log x \cr
0 & \cdots & \cdots & \cdots & 0 & (2\pi i)^n \cr } \in gl_{n+1}(\C).
$$
More precisely,
$$
\Lambda_{j\, k}(x) = \cases{
\ln_k x & when $j=0$ and $k>0$; \cr
{(2\pi i)^j \over (k-j)!}\log^{k-j}x & when $j,k>0$; \cr
0 & when $k < j$.\cr}
$$
We will view this as a multivalued $gl_n(\C)$-valued function on
$\C - \01$. By the {\it principal branch} of $\Lambda(x)$ we shall
mean
the matrix-valued function on the disk $|x - \half|<\half$ obtained
by
taking the standard branches of each of its entries on that disk.
(The
principal branch of $\ln_k$ on this disk is the one given by the
power series expansion (\ref{series}).)

\begin{proposition}\label{diffeq}
The function $\Lambda(x)$ is a fundamental solution of
(\ref{diff-eqn}).
That is,
$$
d\Lambda = \Lambda \bw
$$
and $\Lambda(x)$ is non-singular for each $x\in \C-\01$.\endproof
\end{proposition}

If we analytically continue the principal branch of $\Lambda(x)$
about a loop
in $\C-\01$ based at $\half$, the resulting matrix of functions will
still
be a fundamental solution of (\ref{diff-eqn}). It follows that, for
each loop
$\gamma$ based at $\half$, there is a matrix $M(\gamma)\in
GL_{n+1}(\C)$ such
that the analytic continuation of (the principal branch of)
$\Lambda(x)$
about $\gamma$ is $M(\gamma)\Lambda(x)$.  For a pair of loops
$\alpha,\beta$ based at $\half$, we have
$$
M(\alpha \beta) = M(\alpha) M(\beta).
$$
Since the value of $M(\gamma)$ depends only on the homotopy class of
$\gamma$,
we obtain a monodromy representation
\begin{equation}\label{mono-rep}
M : \pi_1(\C-\01,\half) \to GL_{n+1}(\C).
\end{equation}

Let $\sigma_0,\sigma_1\in \pi_1(\C-\01,\half)$ be the loops defined
by
$$
\sigma_0(t) =  e^{2 \pi i t}/2, \quad \sigma_1(t) = 1 - e^{2\pi i
t}/2,
\quad 0 \le t \le 1.
$$
These loops generate $\pi_1(\C-\01,\half)$.

\begin{proposition}{\rm \cite{ramakrishnan_monod}} \label{monod}
We have
$$
M(\sigma_0) = \left(
\begin{array}{c| c c c}
1 & 0 & \cdots & 0 \\ \hline
0 &   &        &   \\
\vdots & & J & \\
0 & & & \\
\end{array}\right)
\hbox{ and } M(\sigma_1) = \left(
\begin{array}{c| r c c c}
1 & -1 & 0 & \cdots & 0 \\ \hline
0 & 1  &   0  & \cdots  & 0  \\
0 & 0 & 1 &  & 0 \\
\vdots & \vdots &  & \ddots & \\
0 & 0 & 0 & \cdots  & 1\\
\end{array}\right)
$$
where
$$
J = \exp \pmatrix{
0 & 1 & 0 & \cdots & 0 \cr
  & 0 & 1 & \ddots & \vdots \cr
\vdots & & \ddots & \ddots & 0 \cr
  & & & \ddots & 1 \cr
0 & & \cdots & & 0 \cr }
$$
\end{proposition}

\proof The monodromy around $\sigma_0$ is easy to calculate. Indeed,
since
the principal branch of each polylogarithm is single valued on the
disk $|z - 1/2| \le 1/2$, each is unchanged when continued along
$\sigma_0$.
The formula for $M(\sigma_0)$ follows as the analytic continuation of
$\log x$ about $\sigma_0$  is $\log x + 2\pi i$.

Since the principal branch of $\log x$ is defined in a neighbourhood
of
$\R_+$, it follows that $\log x$ is invariant under analytic
continuation
along $\sigma_1$.

We compute the analytic continuation of $\ln_k$ along $\sigma_1$ by
induction. When $k=1$, $\ln_1 x$ changes to $\ln_1 x - 2 \pi i$.
Assume now that $k\ge 1$ and that the continuation of $\ln_k x$
along $\sigma_1$ is
$$
\ln_k x - {2 \pi i \over (k-1)!} \log^{k-1} x.
$$
Denote the integral of $f(z)dz$ along the straight line interval in
the
complex plane  between the points $a$ and $b$ by
$$
\int_a^b f(z) dz.
$$
When  $|x - 1/2| < 1/2$, we have
\begin{eqnarray}\label{first}
\ln_{k+1} x & = & \int_0^x \ln_k z {dz \over z} \nonumber \\
 & = & \int_0^{1/2} \ln_k z {dz \over z} + \int_{1/2}^x \ln_k z {dz
\over z}.
\end{eqnarray}
The result of analytically continuing this along $\sigma_1$ is
\begin{equation}\label{second}
\int_0^{1/2} \ln_k z {dz \over z} + \int_{\sigma_1} \ln_k z {dz \over
z}
+ \int_{1/2}^x \ln_k z {dz \over z}.
\end{equation}
It follows from the inductive formula that the difference between
(\ref{first}) and (\ref{second}) is
\begin{equation}\label{difference}
\int_{\sigma_1} \ln_k z\, {dz \over z} -
{2 \pi i \over (k-1)!} \int_{1/2}^x \log^{k-1} z\, {dz \over z}.
\end{equation}
For each $\epsilon \in (0,1/2)$, the path which traverses the line
segment from $1/2$ to $1-\epsilon$, goes around the boundary of the
disk $|z-1| \le \epsilon$ in the positive direction, then returns
along the interval from $1-\epsilon$ to $1/2$ represents the homotopy
class $\sigma_1$. Again, using the inductive formula for the
monodromy
of $\ln_k x$ around $\sigma_1$, we have
$$
\int_{\sigma_1} \ln_k z\, {dz \over z} =
-{2 \pi i \over (k-1)!}\int_{1-\epsilon}^{1/2} \log^{k-1}z\,
{dz \over z} + \int_{-\pi}^\pi \ln_k(1+\epsilon e^{it})\,
{d(\epsilon e^{it})\over 1 + \epsilon e^{it}}.
$$
The inductive hypothesis also implies that $\ln_k x$ is bounded in a
neighbourhood of 1 when $k>1$, so that the last integral $\to 0$ as
$\epsilon \to 0$ for all $k$. (A separate argument is needed when
$k=1$.)
Combining this with (\ref{difference}), we see that $\ln_{k+1} x$
changes
by
$$
 -{2 \pi i \over (k-1)!}
\left[\lim_{\epsilon \to 0}\int_{1-\epsilon}^{1/2} \log^{k-1}z\,
{dz \over z} + \int_{1/2}^x  \log^{k-1}z\, {dz \over z} \right]
= -{2 \pi i \over k!} \log^k x
$$
when continued around $\sigma_1$. \endproof

The monodromy calculation has several interesting consequences.
First,
even though it is does not make sense, in general, to talk about the
value of a multivalued function at a point, it does make sense to
talk about the value of $\ln_k$ at 1.

\begin{corollary} The value of $(k-1)!\ln_k 1$ is well defined modulo
$\Z(n)$ and is congruent to $(k-1)!\zeta(k)$. \endproof
\end{corollary}

The second important consequence of the monodromy calculation is the
rationality of the monodromy.

\begin{corollary}
The image of the monodromy representation (\ref{mono-rep}) is
contained
in $GL_{n+1}(\Q)$. \endproof
\end{corollary}

The significance of this last result is that it implies that the
local
system over $\C-\01$ which corresponds to the differential equation
(\ref{diff-eqn}) is defined over $\Q$. This local system is called
the
{\it $n$th polylogarithm local system}.  These local systems fit
together
to form an inverse system of local systems whose limit we call the
{\it
polylogarithm local system}. We now describe these local systems in
detail.

Define a meromorphic connection $\nabla$ on the trivial bundle
\begin{equation}\label{bundle}
\P^1 \times \C^{n+1} \to \P^1
\end{equation}
by defining
$$
\nabla f = df - f\bw
$$
where $f : \C - \01 \to \C^{n+1}$ is a section.  This connection has
regular
singular points at $0,1$ and $\infty$, and is flat over $\C-\01$ as
$\bw$
satisfies the integrability condition
$$
d\bw + \bw \wedge \bw = 0
$$
(equivalently, because the equation $\nabla f = 0$ is a system of
ordinary
differential equations).  Let $\bL_0,\bL_1,\ldots,\bL_n$ be
the rows of $\Lambda(x)$. Each of these satisfies (\ref{diff-eqn})
and is
therefore a flat section of (\ref{bundle}).  Even though  these are
multivalued, their $\Q$ linear span is well
defined as the monodromy representation is defined over $\Q$.

Suppose that $X$ is a smooth curve and that $\X$ is a smooth
completion
of $X$. Every flat bundle $E \to X$  has a {\it canonical
extension\/}
$\E \to \X$.   Denote the local monodromy operator about a point
$p \in D := \X - X$ by $T_p$. When each $T_p$ is unipotent, the
canonical
extension is characterized by two properties. First, the meromorphic
extension $\nablabar$ of the connection to $\E \to \X$ has
logarithmic
singularities along $D$. That is,
$$
\nablabar : \overline{\cal E} \to
\overline{\cal E}\otimes_{\O_{\X}} \Omega^1_{\X} (\log D).
$$
Second, the residue of $\nablabar$ at each point of $D$ is nilpotent.

 Denote the $\Q$ local system over $\C-\01$ which corresponds to the
representation (\ref{mono-rep}) by $\V$.

Since $\bw$ has nilpotent residue at 0,1 and $\infty$, we have:

\begin{proposition}\label{poly-loc-sys}
The canonical extension of the flat holomorphic vector bundle
$\V \otimes_\Q \O_{\C -\01}$ to $\P^1$ is the bundle (\ref{bundle})
with the connection $\nabla$ defined above. \endproof
\end{proposition}

\section{The Bloch-Wigner function}
\label{bw_function}

Define
$$
D_2(x) = \Im \ln_2 x + \log |x| \arg(1-x)
$$
when $|x-1/2| < 1/2$ and where $\ln_2 x$, $\log x$, and $\arg (1-x)$
denote the principal branches of these functions in the disk
$|x-1/2| < 1/2$.  An easy computation using the monodromy calculation
(\ref{monod}) shows that $D_2$ is invariant under continuation along
the generators $\sigma_0$ and $\sigma_1$ of $\pi_1(\C-\01,1/2)$.
Consequently, the function $D_2$ extends to a single valued, real
analytic
function
$$
D_2 : \C - \01 \to \R.
$$
This is called the {\it Bloch-Wigner function\/}. If we define
$$
D_2(0) = D_2(1) = D_2(\infty) = 0
$$
then $D_2$ is a continuous function $D_2 : \P^1 \to \R$.  The Bloch-%
Wigner function should be viewed as having the same relation to
$\ln_2$
as $D_1 := \log |\phantom{x}|$ bears to the logarithm.

The boundary of hyperbolic 3-space $\H^3$ is the Riemann sphere
$\P^1$.
The group of orientation preserving isometries of $\H^3$ is
$PSL_2(\C)$.
The induced action on the boundary is just the standard action of
$PSL_2(\C)$ on $\P^1$ via fractional linear transformations.

Denote the ideal tetrahedron in $\H^3$ with vertices at
$a_0, a_1, a_2, a_3 \in \P^1$ by $\Delta(a_0,a_1,a_2,a_3)$. Since
the volume form of hyperbolic space is invariant under the action of
the isometry group,
$$
\vol \Delta(a_0,a_1,a_2,a_3) = \vol \Delta (\lambda,1,0,\infty).
$$
where $\lambda$ is the cross ratio $[a_0:a_1:a_2:a_3]$ of the
vertices.

The following result goes back to Lobachevsky (cf. \cite{milnor}).
A proof may be found in \cite[p.~172]{dupont-sah}.

\begin{theorem}\label{volume} For each $z \in \P^1$,
the volume of $\Delta (z,1,0,\infty)$ equals $D_2(z)$.
\end{theorem}

\begin{corollary}\label{funct_eqn_1}
If $a_0, a_1, a_2, a_3, a_4 \in \P^1$, then
$$
\sum_{j=0}^4 (-1)^j D_2([a_0:\cdots:\widehat{a_j}:\cdots:a_4]) = 0.
$$
Moreover, for all permutations $\sigma$ of $\{0,1,2,3\}$,
$$
D_2([a_{\sigma(0)}:\cdots:a_{\sigma(3)}]) =
\sgn(\sigma)D_2([a_0:\cdots:a_3]).
$$
\end{corollary}

\proof
The first assertion follows as the ideal polyhedron $P$ with vertices
$a_0, a_1, a_2, a_3, a_4$ decomposes into a union of ideal tetrahedra
in 2 different ways. Viz.
$$
P = \Delta(a_1, a_2, a_3, a_4) \cup \Delta(a_0, a_1, a_3, a_4)
\cup \Delta(a_0, a_1, a_2, a_3)
$$
and
$$
P = \Delta(a_0, a_2, a_3, a_4) \cup \Delta(a_0, a_1, a_2, a_4)
$$
In each case the pieces intersect along 2-dimensional faces. The
first
assertion follows from  Theorem \ref{volume} by comparing volumes.
The second follows as swapping the order of 2 vertices reverses the
orientation of the tetrahedron. \endproof
\medskip

Taking the five points to be $y,x,1,0,\infty$, we obtain the usual
form
of the functional equation, which is the analogue for $D_2$ of the
Abel-Spence functional equation of $\ln_2$.

\begin{corollary}\label{funct_eqn_2}
If $y,x, 1, 0, \infty$ are distinct points of $\P^1$, then
$$
D_2(x) - D_2(y) + D_2(y/x) - D_2((1-y)/(1-x)) +
D_2((1-y^{-1})/(1-x^{-1})) = 0. \quad \square
$$
\end{corollary}

Ramakrishnan \cite{ramakrishnan_bw} showed that all the
polylogarithms
have such single valued cousins. The essential point being that one
can use the
unipotence of the monodromy group and induction to kill off the
monodromy
of $\ln_k$. Zagier \cite{zagier} gave an explicit
formula for Ramakrishnan's functions. In general, there seems to be
no
canonical way to go from the multivalued polylogarithm to a single
valued
function.  Goncharov \cite{goncharov} has modified Ramakrishnan's
trilogarithm and proved that his function satisfies a very natural
functional equation. His function is defined by
\begin{equation}\label{gonch_trilog}
D_3(x) = \Re\left[ \ln_3(x) - \log|x| \ln_2 x + (\log^2|x| \ln_1
x)/3\right].
\end{equation}
It is single valued and continuous on $\P^1$, and real analytic on
$\P^1 - \{0,1,\infty\}$.

\section{The regulator $K_3(\protect\C) \to
\protect\C/\protect\Q(2)$}
\label{reg-k3}
The Deligne cohomology of $\spec \C$ is
$$
\Hdel^m(\spec \C,\Lambda(k)) = \cases{
0 & $m \neq 1$; \cr
\C/\Lambda(k) & $m = 1$, \cr}
$$
where $\Lambda$ denotes $\Z$, $\Q$ or $\R$. So the only non-trivial
regulators for $\spec \C$ with values in Deligne cohomology are
$$
K_{2m-1}(\C) \to \Hdel^1(\spec \C, \Lambda(m)) \approx \C/\Lambda(m)
$$
for each $m \ge 1$. The first regulators
$$
K_1(\C) \approx \C^\ast \to \C/\Z(1)\quad\hbox{and}\quad
K_1(\C) \approx \C^\ast \to \C/\R(1) \approx \R
$$
are given by the functions $\log$ and $\log|\phantom{x}|$,
respectively.
In this section we give an account of the construction of the second
regulators
$$
K_3(\C) \to \C/\Q(2) \quad \hbox{and}
\quad K_3(\C) \to \C/\R(2) \approx \R(1)
$$
using $\ln_2$ and $D_2$ respectively. These results go back to Bloch
and
Wigner (unpublished, cf. \cite{dupont-sah}). As is customary,
$GL_n(\C)^\delta$ signifies that $GL_n(\C)$ is viewed as a group with
the discrete topology.

\begin{lemma}\label{gen_posn}
The Bloch-Wigner function defines a canonical group cohomology class
$$
\D_2 \in H^3(GL_2(\C)^\delta,\R).
$$
\end{lemma}

\proof Let $F$ be a field and $n\in \N$. For each $k\in\N$, define
$C_k(F,n)$ to be the free abelian group with basis ordered
$(k+1)$-tuples $(v_0,\ldots,v_k)$, of vectors $v_j$ of $F^n$ in
general
position. (I.e. , each $\min(n,k+1)$ of them are independent.) Define
a differential $\partial : C_k \to C_{k-1}$ by
$$
\partial : (v_0,\ldots,v_k) \mapsto
\sum_{j=0}^k (-1)^j (v_0,\ldots,\widehat{v_j},\ldots,v_k).
$$
When $F$ is infinite, the complex
$C_\bullet(F,n)$ is quasi-isomorphic to $\Z$. Since $GL_n(F)$ acts on
this complex, it is a resolution of the trivial module. So, for all
$GL_n(F)$ modules $M$, there is a natural map
$$
H^\bullet\left(\Hom_{GL_n(F)}(C_\bullet(F,n),M)\right)
\to H^\bullet(GL_n(F),M),
$$
provided that $F$ is infinite.

Denote the point in $\P^{n-1}(F)$ determined by the non-zero vector
$v$ of $F^n$ by $[v]$.

Define a map $f : C_3(\C,2) \to \R $ by
$$
f(v_0,v_1,v_2,v_3) =
D_2\left(\left[[v_0]:[v_1]:[v_2]:[v_3]\right]\right).
$$
The cocycle condition is simply the functional equation
(\ref{funct_eqn_1}).
Denote the image of this cohomology class in
$H^3(GL_2(\C)^\delta,\R)$
by $\D_2$. \endproof
\medskip

A cohomology class $c \in H^m(GL_m(R),\Lambda)$ defines a map
$$
K_m( R) \to \Lambda.
$$
This map is obtained as the composite
$$
K_m(R) = \pi_m(BGL(R)^+) \to H_m(BGL(R)^+) \approx H_m(GL(R))
\mapright{c}\Lambda
$$
of the Hurewicz homomorphism with the $\Lambda$-valued functional
on homology induced by $c$.

\begin{theorem}\label{real_reg}
The map
$$
K_3(\C) \to \C/\R(2)
$$
defined by $i\D_2$ equals the Beilinson regulator, and this map
equals
half of the Borel regulator.
\end{theorem}

\proof The proof is an assemblage of results from the literature, and
we only give a sketch. First, the cocycle
$$
GL_2(\C)^4 \to \C/\R(2)
$$
induced by $iD_2$, composed with the cross ratio, represents a
continuous
cohomology class $i\gD_2$, even though this cocycle is not itself
continuous.  This apparently contradictory state of affairs arises
because
the cross ratio $\left(\P^1\right)^4 \to \P^1$ is not everywhere
defined.
One can show (e.g. \cite{yang_thesis}) that the image of this class
under the natural map
$$
\Hcts^3(GL_2(\C),\C/\R(2)) \to H^3(GL_2(\C)^\delta,\C/\R(2))
$$
is the class $i\D_2$ of (\ref{gen_posn}).

Next, since
$$
\Hcts^3(GL_2(\C),\C/\R(2)) \approx \Hcts^3(GL_n(\C),\C/\R(2))
\approx \C/\R(2)
$$
for all $n \ge 2$, and since this group is spanned by the continuous
cohomology class $\beta_2$ used to define the Borel regulator (cf.
\cite{ramakrishnan_survey}, \cite{D-H-Z}), it follows that
there is a real number $\lambda$ such that
$$
i \gD_2 = \lambda \beta_2/2.
$$

Denote by $c_k$ the $k$th Beilinson Chern class of the universal flat
bundle over $BGL_n(\C)^\delta$. This is an element of
$$
\Hdel^{2k}(BGL_n(\C)^\delta,\R(k)) \approx
H^{2k-1}(GL_n(\C),\C/\R(k)).
$$
In \cite{D-H-Z} it is shown that, for all $n$ and $k$, the classes
$c_k$ and $\beta_k/2$ in
$$
H^k(GL_n(\C)^\delta,\C/\R(2))
$$
are equal. It follows that $i \gD_2$ is $\lambda$ times the
Beilinson Chern class
$$
c_2 : K_3(\C) \to \C/\R(2).
$$

Denote the rational $K$-theory $K_\bullet(R)\otimes \Q$ of a ring $R$
by $K_\bullet(R)_\Q$.
Since $BGL(R)^+$ is an H-space, the Hurewicz homomorphism
$$
K_m(R)_\Q := \pi_m(BGL(R)^+)\otimes \Q \to
H_m(BGL(R)^+,\Q) \approx H_m(GL(R),\Q)
$$
is injective. Define the {\it rank filtration}
$$
r_1K_m(R) \subseteq r_2K_m(R) \subseteq r_3K_m(R) \subseteq \cdots
\subseteq K_m(R)_\Q
$$
of $K_m(R)_\Q$ by
$$
r_kK_m(R) = K_m(R)_\Q \cap
\im\left\{H_m(GL_k(R),\Q) \to H_m(GL(R))\right\}.
$$
Suslin \cite{suslin} has proved that for all infinite fields $F$
$$
K_m(F)_\Q = r_m K_m(F)
$$
and that
$$
K_m(F)_\Q/r_{m-1}K_m(F) \approx K_m^M(F)\otimes \Q,
$$
where $K_\bullet^M(F)$ denotes the Milnor $K$-theory of $F$.

In particular, for all infinite fields $F$, there is a canonical
isomorphism
$$
K_3(F)_\Q = K_3^M(F)_\Q \oplus r_2K_3(F).
$$
Since all elements of $K_3^M$ of any field are decomposable, the
generalization of the Whitney sum formula to $K$-theory implies
that
$$
c_2 : K_3(\C) \to \C/\R(2)
$$
vanishes on $K_3^M(\C)$. It follows that $c_2$ factors through
the projection onto $r_2K_3(\C)$.
$$
\matrix{
K_3(\C) & \mapright{c_2} & \C/\R(2) \cr
\mapdownalt{{\rm proj}} & \nearrow & \cr
r_2K_3(\C) & & \cr}
$$
Since $r_2K_3(\C)$ comes from $H_3(GL_2(\C))$, the restriction of
$c_2$ to $r_2K_3(\C)$ is given by the class $i\lambda^{-1}\gD_2$.
Now Dupont \cite{dupont} has proved that the $i\gD_2$ equals the
Cheeger-Simons class
$$
\hat{c}_2 \in H^3(GL_2(\C)^\delta,\C/\R(2)).
$$
of the universal rank 2 flat bundle. But by the main result of
\cite{D-H-Z}, this class equals the Beilinson Chern class of the
universal flat bundle. It follows that $\lambda = 1$. \endproof

Most of this story has been extended to the trilogarithm by
Goncharov \cite{goncharov} and Yang \cite{yang_thesis}. In
\cite{hain-macpherson} a canonical single valued trilogarithm
$$
S_3 : \left\{ \matrix
{\hbox{ ordered 6-tuples of points \hfil} \cr
\hbox{ in $\P^2$, no 3 on a line \hfil} \cr} \right\}/
\hbox{projective equivalence} \to \R
$$
was constructed. It satisfies  the seven term functional equation
$$
\sum_{j=0}^6 (-1)^j S_3(a_0,\ldots,\widehat{a_j},\ldots,a_6) = 0
$$
where $a_0,\ldots,a_6$ are points in $\P^2$, no 3 of which lie on a
line. This equation is an obvious generalization of the 5-term
equation satisfied by $D_2$.  As in (\ref{gen_posn}), this function
determines a class in $H^5(GL_3(\C),\R)$.

As in the case of the dilogarithm, the cocycle condition is precisely
the functional equation.  Goncharov and Yang both
showed that this class in $H^5(GL_3(\C),\R)$ is a non-zero rational
multiple of the class which corresponds to the Beilinson Chern class
of the universal flat bundle over $BGL_3(\C)^\delta$. Goncharov
found an explicit formula for $S_3$ in terms of his single valued
version of the classical trilogarithm (\ref{gonch_trilog}). The
appropriate analogue of Suslin's theorem is not known for $K_5(\C)$.
However, Yang \cite{yang} proved the rank conjecture
for all $K$-groups of all number fields except $\Q$.  Their
work enables one to find formulas for the regulator mapping
$$
c_3 : K_5(F) \to \left[\C/\R(3)\right]^{r_1 + r_2}
$$
for all number fields $F$ in terms of values of the trilogarithm of
\cite{hain-macpherson} in the case of Yang, and the classical
trilogarithm in the case of Goncharov.

The regulator $c_2 : K_3(\C) \to \C/\Q(2)$ can be written in
terms of the multivalued dilogarithm. This work goes back to
unpublished
work of Bloch and Wigner (cf. \cite{dupont-sah}). The construction is
similar to, but more complicated than, the construction of the
regulator
given above. We only sketch the construction. More details can be
found
in \cite{dupont-sah}.

For $x$ in the disk $|x-1/2|< 1/2$, define
$$
\rho(x) = {1\over 2}\left[\log x \wedge \log(1-x) +
2\pi i \wedge {1 \over 2 \pi i}
\left( \ln_2(1-x) -\ln_2(x) - \pi^2/6\right)\right]
\in \Lambda^2_\Z \C.
$$
Here all functions are taken to be the principal branches. This
function
is single valued, and therefore extends to a single valued function
$$
\rho : \C-\01 \to \Lambda^2_\Z \C.
$$
It satisfies a generalization of the 5-term equation satisfied
by $D_2$. If $x,y \in \C- \01$ and $x \neq y$, then
$$
\rho(x) - \rho(y) + \rho(y/x) - \rho((1-y)/(1-x)) +
\rho((1-y^{-1})/(1-x^{-1})) = 0.
$$
Define $\cP(F)$ to be the free abelian group generated
by $F -\01$ subject to the relations
$$
[x] - [y] + [y/x] - [(1-y)/(1-x)] + [(1-y^{-1})/(1-x^{-1})] = 0.
$$
This is often called the {\it scissors congruence group}, or the {\it
Bloch
group}.  The function $\rho$ induces a map
$$
\rho : \cP(\C) \to \Lambda^2_\Z\C.
$$

There is a natural homomorphism
$$
H_3(SL_2(F)) \to \cP(F)
$$
whose construction is similar to that of the homomorphism
$H_3(GL_2(F)) \to H_\bullet(C_\bullet(F,n))$
given in the proof of (\ref{gen_posn}).
If $F$ is algebraically closed of characteristic 0, then there is
an exact sequence
$$
0 \to \Q/\Z \to H_3(SL_2(F)) \to \cP(F) \to \Lambda^2_\Z F^\times
\to K_2(F) \to 0.
$$
A proof can be found in the appendix of \cite{dupont-sah}. The map
$\cP(F) \to \Lambda^2_\Z F^\times$ takes the generator $[x]$ of
$\cP(F)$ to
$(1-x)\wedge x$. The most right hand map takes $x\wedge y$ to
$\{x,y\}$.

The kernel of the map
$$
\Lambda^2 \exp : \Lambda_\Z^2 \C \to \Lambda_\Z^2 \C^\ast
$$
is $\C/\Q(2)$, where this is included in $\Lambda_\Z^2 \C$ by taking
the
coset of $\lambda$ to $2\pi i \wedge (\lambda/2\pi i)$. Since the
diagram
$$
\matrix{
\cP(\C) & \to & \Lambda^2_\Z \C^\times \cr
\mapdownalt{\rho}& & \mapdown{id} \cr
\Lambda_\Z^2 \C & \mapright{\wedge^2\exp} & \Lambda_\Z^2 \C^\ast}
$$
commutes, $\rho$ induces a homomorphism
$$
H_3(SL_2(\C)) \to \C/\Q(2).
$$
By \cite{dupont}, this represents the second Cheeger-Simons class
$c_2$
of the universal flat bundle over $BSL_2(\C)^\delta$. By the main
result
of \cite{D-H-Z}, this equals the Beilinson Chern class of the
universal
flat bundle over $BSL_2(\C)^\delta$.  As above, $c_2$ vanishes on
$K_3^M(\C)$ and the diagram
$$
\matrix{
K_3(\C) & \mapright{c_2} & \C/\Q(2) \cr
\mapdownalt{{\rm proj}} & \mapne{\rho} & \cr
r_2K_3(\C) & & \cr}
$$
commutes. One can show that the map
$$
H_3(SL_2(F),\Q) \to r_2K_3(F)
$$
is surjective. It follows that the map constructed above induces the
regulator on all of $K_3(\C)$. Alternatively, one can appeal to the
theorem of Suslin \cite{suslin_K3ind} which asserts that there are
natural isomorphisms
$$
H_3(SL_2(F),\Q) \approx r_2K_3(F)\approx
\ker \left\{\cP(F) \to \Lambda^2_\Z F^\times \right\}\otimes \Q
$$
for all fields $F$. This last isomorphism says that all of the weight
2
part of $K_3$ of a field comes from $\P^1 - \{0,1,\infty\}$ and that
all the relations come from the functional equation of the
dilogarithm.
This generalizes the fact mentioned in the introduction that the
relations
in $K_1$ come from the functional equation of the logarithm. The
analogue of this statement for the weight 3 part of $K_5$ is not
known
at this time, although Goncharov \cite{goncharov} has made
significant
progress.

\section{Iterated integrals}
\label{it-ints}

At this stage it is convenient to introduce Chen's iterated integrals
\cite{chen} Basic references for this section are \cite{chen},
\cite{geom}.

Suppose that $M$ is a manifold and that $\row w 1 r$ are smooth
$\C$-valued
1-forms on $M$.
For each piecewise smooth path $\gamma : [0,1] \to M$, we can define
$$
\int_\gamma \row w 1 r =
\int \cdots \int_{0 \le t_1 \le \cdots \le t_r \le 1}
f_1(t_1) \cdots f_r(t_r)\, dt_1 \ldots dt_r,
$$
where $\gamma^\ast w_j = f_j(t)\, dt$ for each $j$.
This can be viewed as a $\C$-valued function
$$
\int_\gamma \row w 1 r : PM \to \C
$$
on the the space of piecewise smooth paths in $M$.
When $r=1$, $\int_\gamma w$ is just the usual line integral.
An {\it iterated integral} is any function $PM \to \C$ which is a
linear
combination of constant functions and basic iterated integrals
$$
\int_\gamma \row w 1 r.
$$

Now let $M = \C - \01$ and
$$
\w_0 = {dz \over z}\hbox{ and } \w_1 = {dz \over 1-z}.
$$
Then
$$
\ln_1 x = - \log(1-x) = \int_0^x \w_1.
$$
By induction and the definition, we have, for all $k\ge 2$,
$$
\ln_k x = \int_0^x \ln_{k-1} z\, \w_0 =
\int_0^x \w_1\overbrace{\w_0 \ldots \w_0}^{k-1\rm \; times}.
$$
Here the path of integration must be chosen so that once it has left
0,
it never passes through 0 or 1 on its way to $x$.

The basic properties of iterated integrals are summarized in the
following proposition.

\begin{proposition}{\rm \cite{chen,geom}}
\label{props}
Suppose that $w_1, w_2, \ldots $ are $\C$-valued 1-forms on a
manifold $M$.
\begin{enumerate}
\item[(i)] The value of $\int_\gamma \row w 1 r$ is independent of
the
parameterization of $\gamma$.
\item[(ii)] If $\alpha,\beta : [0,1] \to M$ are composable paths
(i.e. $\alpha(1) = \beta(0)$), then
$$
\int_{\alpha \beta} \row w 1 r =
\sum_{j=0}^r \int_\alpha \row w 1 i \int_\beta \row w {i+1} r.
$$
Here, $\int_\gamma \row \phi 1 m$ is to be interpreted as 1 when
$m=0$.
\item[(iii)] For all paths $\gamma$,
$$
\int_{\gamma^{-1}} \row w 1 r = (-1)^r \int_\gamma \row w r 1.
$$
\item[(iv)] For all paths $\alpha$ in $M$,
$$
\int_\alpha \row w 1 r \int_\alpha \row w {r+1} {r+s} =
\sum_\sigma \int_\alpha \row w {\sigma(1)} {\sigma(r+s)},
$$
where $\sigma$ ranges over all shuffles of type $(r,s)$.
\end{enumerate}
\end{proposition}

\section{The Regulator $K_2(X) \to H^1(X,\protect\C^\ast)$}
\label{heisenberg}

This section is an exposition of Bloch's construction of the
regulator
$$
c_2 : K_2(X) \to \Hdel^2(X,\Z(2))
$$
using the dilogarithm. We have freely incorporated the elegant
approaches of Deligne \cite{deligne_letter} and
Ramakrishnan \cite{ramakrishnan_heisenberg,ramakrishnan_survey}.

When $X$ is a curve, there is a natural isomorphism
$$
H^2(X,\Z(2)) \approx H^1(X,\C/\Z(2)).
$$
Identifying $\C/\Z(2)$ with $\C^\ast$ by the map
$\lambda \mapsto \exp[\lambda/2\pi i]$, we obtain a canonical
identification of $\Hdel^2(X,\Z(2))$ with $H^1(X,\C^\ast)$, the group
of flat line bundles over $X$.

We set
$$
H_\Z = \pmatrix{
1 & \Z(1) & \Z(2) \cr
0 & 1     & \Z(1) \cr
0 & 0 & 1 \cr }
$$
and
$$
H_\C= \pmatrix{
1 & \C & \C \cr
0 & 1     & \C \cr
0 & 0 & 1 \cr }
$$
There is a natural bundle projection
\begin{equation}\label{heis-bundle}
H_\Z\backslash H_\C \to \C^\ast \times \C^\ast;
\end{equation}
it takes the coset of
$$
\pmatrix{
1 & u & w \cr
0 & 1 & v \cr
0 & 0 & 1 \cr}
$$
to $(e^u,e^v)$. It has fiber $\C/\Z(2)$, which we identify with
$\C^\ast$
as above.

\begin{proposition}\label{connection}
There is a natural connection on this bundle with curvature
$(dx/x) \wedge (dy/y)/2 \pi i$, where $x$ and $y$ are the coordinates
in
$\C^\ast \times \C^\ast$.
\end{proposition}

\proof First consider the pullback of the bundle (\ref{heis-bundle})
to $\C\times\C$
\begin{equation}\label{lift}
\Z(2)\backslash H_\C \to \C \times \C
\end{equation}
along the map $(u,v) \mapsto (e^u,e^v)$.  The map
$$
H_\C \to \C^\ast \times \C \times \C
$$
defined by
$$
\pmatrix{
1 & u & w \cr
0 & 1 & v \cr
0 & 0 & 1 \cr}
\mapsto (\exp(w/2\pi i), u,v)
$$
induces an isomorphism of $\Z(2)\backslash H_\C$ with
$\C^\ast \times \C \times \C$ which commutes with the projections to
$\C\times \C$. So the bundle (\ref{lift}) is trivial, and sections
of it can be identified with maps $\zeta : \C \times \C \to \C^\ast$.
Define a connection on this bundle by
$$
\nabla \zeta = d \zeta - \zeta udv/2 \pi i
$$
This connection is easily seen to be invariant under the left action
$$
(n,m) : (\zeta, u,v) \mapsto (e^{nv}\zeta, u + 2\pi i n, v + 2\pi i
m)
$$
of $\Z\times \Z$ on $\C^\ast \times \C \times \C$ induced by the left
action of $H_\Z$ on $H_\C$. It therefore descends to a connection on
the bundle (\ref{heis-bundle}).  The connection form of the pullback
bundle is $udv/2\pi i$, from which it follows that its curvature is
$du\wedge dv/2\pi i$ and that the curvature of (\ref{heis-bundle}) is
$(dx/x) \wedge (dy/y)/2\pi i$. \endproof

Now suppose that $X$ is a smooth curve over $\C$. Denote the function
field of $X$ by $\C(X)$ and the generic point $\spec \C(X)$ of $X$ by
$\eta_X$. We first define the regulator on $K_2(\eta_X) :=
K_2(\C(X))$.
The Deligne cohomology of $\eta_X$ is defined by
$$
\Hdel^m(\eta_X,\Lambda(k)) = \lim_{\to} \Hdel^m(U,\Lambda(k))
$$
where the limit is taken over all Zariski open subsets $U$ of $X$. In
particular, we have
$$
\Hdel^2(\eta_X,\Z(2)) = \lim_{\to} H^1(U,\C^\ast).
$$
This latter group is the group of flat line bundles at the generic
point
of $X$---elements of this group are flat line bundles defined on some
Zariski open subset of $X$, and two such are identified if they agree
on
a smaller open subset. The product is tensor product.

By Matsumoto's Theorem \cite{milnor_book}, $K_2(\eta_X)$ is generated
by symbols $\{f,g\}$, where $f,g \in \C(X)^\times$. The only
relations
which hold between these symbols are bilinearity
$$
\{f_1f_2,g\} = \{f_1,g\}\{f_2,g\},\quad \
\{f,g_1g_2\} = \{f,g_1\}\{f,g_2\}
$$
and the Steinberg relation
$$
\{1-f,f\} = 1
$$
whenever $f$ and $1-f$ are both in $\C(X)^\times$.

Now suppose that $\{f,g\}\in K_2(\eta_X)$. There is a Zariski open
subset $U$ of $X$ such that $f$ and $g$ are both defined and
invertible
on $U$.  They therefore define a regular function
$$
(f,g) : U \to \C^\ast \times \C^\ast.
$$
The pullback of the line bundle $H_\Z\backslash H_\C$ to $U$ is flat
as it has curvature a multiple of $(df/f)\wedge (dg/g)$, which is
zero
as $U$ is a curve.
Denote it by $\langle f,g\rangle$.  It is an element of
$H^1(U,\C^\ast)$,
and therefore of $H^1(\eta_X,\C^\ast)$.

\begin{remark}\label{natural}
Observe that this construction makes sense when $X$ is a Riemann
surface and $\C(X)$ denotes the field of meromorphic functions of
$X$.
It is natural with respect to holomorphic maps between Riemann
surfaces.
\end{remark}

\begin{proposition}\label{cocycle}
If $f,g$ are invertible functions on the Zariski open subset $U$ of
$X$, the monodromy of the flat bundle $\langle f,g \rangle$ about
a loop $\gamma$ based at $p\in U$ is
$$
I(f,g,\gamma) :=  \int_\gamma \dlog f \dlog g -
\log g(p) \int_\gamma \dlog f + \log f(p) \int_\gamma \dlog g
\in \C/\Z(2).
$$
\end{proposition}

\proof We will deduce the result by computing the monodromy about
a loop $\gamma$ in $\C^\ast \times \C^\ast$. The assertion will then
follow by pulling back the answer to $U$. Let $\gamma$ be a path in
$\C^\ast \times \C^\times$ which begins at $(x_0,y_0)$.
For $t \in [0,1]$, denote the path $s \mapsto \gamma(st)$ by
$\gamma_t$. The horizontal lift of $\gamma$ to $H_\Z\backslash H_\C$
which begins at the coset of
$$
\pmatrix{
1 &  \log x_0 & w \cr
0 & 1 & \log y_0 \cr
0 & 0 & 1 }
$$
is
$$
t \mapsto H_\Z
\pmatrix{
1 &  \log x_0 & w \cr
0 & 1 & \log y_0 \cr
0 & 0 & 1 }
\pmatrix{
1 &  \int_{\gamma_t} \dlog x &  \int_{\gamma_t} \dlog x \dlog y \cr
0 & 1 & \int_{\gamma_t} \dlog y \cr
0 & 0 & 1 \cr}.
$$
That is,
$$
t \mapsto \pmatrix{
1 & \log x_0 + \int_{\gamma_t} \dlog x & w +
\int_{\gamma_t} \dlog x \dlog y + \log x_0 \int_{\gamma_t} \dlog y
\cr
0 & 1 & \log y_0 + \int_{\gamma_t} \dlog y \cr
0 & 0 & 1 \cr}
$$
Now suppose that $\gamma$ is a loop. Then the endpoint of the
horizontal
lift of $\gamma$ is congruent to the matrix
$$
\pmatrix{
1 & -\int_\gamma \dlog x & \int_{\gamma} \dlog x \int_{\gamma} \dlog
y \cr
0 & 1 & -\int_\gamma \dlog y \cr
0 & 0 & 1 \cr}
\pmatrix{
1 & \log x_0 + \int_{\gamma} \dlog x & w +
\int_{\gamma} \dlog x \dlog y + \log x_0 \int_{\gamma} \dlog y \cr
0 & 1 & \log y_0 + \int_{\gamma} \dlog y \cr
0 & 0 & 1 \cr}
$$
$$
=
\pmatrix{
1 & \log x_0 & w + \int_{\gamma} \dlog x \dlog y +
\log x_0 \int_{\gamma} \dlog y - \log y_0 \int_\gamma \dlog x \cr
0 & 1 & \log y_0 \cr
0 & 0 & 1 \cr}
$$
modulo the left action of $H_\Z$. It follows that the holonomy about
$\gamma$ is
$$
 \int_{\gamma} \dlog x \dlog y + \log x_0 \int_{\gamma} \dlog y -
\log y_0 \int_\gamma \dlog x \hbox{ mod }{\Z(2)}.
$$
The result follows by pulling the result back to $U$ along the map
$(f,g)$. \endproof

\begin{proposition}\label{symbol}
If $f,f_1,f_2,g \in \C(X)^\times$, then
$$
\langle f_1 f_2, g \rangle = \langle f_1,g \rangle \langle f_2, g
\rangle
\quad \hbox{and}\quad \langle f,g \rangle^{-1} = \langle g,f \rangle.
$$
Moreover, if $f,1-f\in \C(X)^\ast$, then
$$
\langle 1-f,f \rangle = 1.
$$
\end{proposition}

\proof
It is clear from properties of the logarithm and (\ref{props}) that
$$
I(f_1f_2,g, \gamma) = I(f_1,g,\gamma) + I(f_2,g,\gamma).
$$
It is not difficult to use (\ref{props}) to prove that
$$
I(f,g,\gamma) + I(g,f,\gamma) = 0.
$$
These imply the linearity and skew symmetry of the symbol
$\langle f,g \rangle$.

It suffices to prove the Steinberg relation in the universal case
where $U = \C - \01$ and $f=x$. The  line bundle
$\langle 1-x,x \rangle$ is trivial as a flat bundle if and only if it
has a flat section. The dilogarithm provides such a section.
Define $s : \C - \01 \to H_\Z\backslash H_\C$ by
$$
s(x) = H_\Z\pmatrix{
1 & \log(1-x) & - \ln_2 x \cr
0 & 1 & \log x \cr
0 & 0 & 1 \cr}
$$
This section is flat, so the Steinberg relation holds. \endproof

\begin{theorem}{\rm{\cite{beilinson}}}
Taking $\{f,g\}$ to $\langle f,g \rangle$ defines a map
$$
K_2(\eta_X) \to H^1(\eta_X,\C^\ast)
$$
which is the Chern class $c_2$.
\end{theorem}

\proof The first assertion is an easy consequence of Matsumoto's
description of $K_2$ and (\ref{symbol}). The second follows from the
fact
that the symbol $\{f,g\}$ is the cup product of
$f,g\in K_1(\eta_X)\approx \C(X)^\times$.
Properties of Chern classes then imply that
$$
c_2(\{f,g\}) = c_1(f) \cup c_1(g)
$$
where the right hand side is the cup product of
$$
c_1(f), c_2(g) \in \Hdel^2(\eta_X,\Z(1)) \approx \C(X)^\ast.
$$
Under this isomorphism, $c_1$ is just the identity.

The formula for the cup product in Deligne cohomology implies that
$c_1(f)\cup c_1(g)$ is represented by the element of
$$
\Hdel^2(\eta_X,\Z(2)) \approx H^1(\eta_X,\C/\Z(2))
$$
defined by
$$
\gamma \mapsto I(f,g,\gamma).\quad \square
$$

We now globalize this construction. For each $x\in X$, there is map
$$
\delta_x : H^1(\eta_X,\C^\ast) \to \C^\ast.
$$
To define $\delta_x(l)$, represent $l$ by  a flat line bundle $L \to
U$
over a Zariski open subset $U$ of $X$. Chose a small closed disk
$\Deltabar$ in $X$, centered at $x$, such that $\Deltabar - \{x\}$ is
contained in $U$. Define $\delta(l)$ to be the monodromy of $L$ about
the boundary of $\Deltabar$.

Suppose that $\nu : F^\times \to \Z$ is a valuation on a field $F$.
Let $\cal O$ be the associated valuation ring (i.e., 0 and those
elements
of $F^\times$ with valuation $\ge 0$). Let $\gP$ be the maximal
ideal of $\cal O$.  The tame symbol of $f,g \in F^\times$ is defined
by
$$
(f,g)_\nu =
(-1)^{\nu(f)\nu(g)} {f^{\nu(g)}\over g^{\nu(f)}}\hbox{ mod } \gP
$$
(See, e.g., \cite{milnor}.)
To each $x\in X$ associate the valuation which takes a function $f$
to its order $\nu_x(f)$ at $x$. In this way we associate a tame
symbol
$(\phantom{x},\phantom{x})_x$ to each $x\in X$.

\begin{proposition}\label{tame}
Suppose that $x \in X$. If $f,g \in \C(X)$, then
$$
\delta_x\langle f,g \rangle = (f,g)_x.
$$
\end{proposition}

\proof Since both sides of the expression in the statement of the
proposition are skew symmetric and bilinear, we can reduce, with the
help of (\ref{natural}),
to the following 2 cases. First, if $z$ is a local holomorphic
parameter
about $x$, then we have to show that $\delta_x\langle z,z \rangle =
-1$.
Second, if $f$ is a unit in a neighbourhood of $x$, then
$$
\delta_x\langle f,g \rangle = f^{\nu_x(g)}(x)
$$

{}From (\ref{cocycle}) and (\ref{props})(ii) it follows that
$$
I(z,z,p) = \int_{\partial \Delta} \dlog z \dlog z =  {1\over 2}
\left(\int_{\partial \Delta}\dlog z\right)^2
= {(2\pi i )^2 \over 2} \in \C/\Z(2),
$$
where $\Delta$ is a sufficiently small imbedded disk in $X$ centered
at $x$.
Under the standard isomorphism $\C/\Z(2) \approx \C^\ast$, this
corresponds to $e^{i\pi}= -1$. This proves the first assertion.

To prove the second, write $f = e^\phi$, where $\phi$ is holomorphic
in
a neighbourhood of $x$. By (\ref{cocycle}), we have
$$
I(f,g,p) =  \int_{\partial \Delta} (\phi(z) - \phi(p))\dlog g
+ \phi(p)\int_{\partial \Delta}\dlog g =
  \int_{\partial \Delta}\phi(z)\dlog g.
$$
By the Residue Theorem this equals $2 \pi i \nu_x(g)\phi(x)$.
So
$$
\delta_x \langle f,g \rangle = \exp\left(\nu_x(g)\phi(x)\right)=
f(x)^{\nu_x(g)}.
$$
This proves the second assertion. \endproof

The following version of the Gysin sequence
is easily verified.

\begin{proposition}\label{gysin}
The sequence
$$
0\to H^1(X,\C^\ast) \to H^1(\eta_X,\C^\ast)
\mapright{\oplus \delta_x} \bigoplus_{x\in X}\C\ast
$$
is exact. \endproof
\end{proposition}

The analogue of the Gysin sequence in algebraic $K$-theory is the
localization sequence (reference !). In our case, it asserts that
the sequence
$$
K_2(X) \to K_2(\eta_X)
\mapright{\oplus (\phantom{x},\phantom{x})_x}
\bigoplus_{x\in X} \C^\ast
$$
is exact.  By (\ref{tame}), the diagram
$$
\matrix{ & &
K_2(X) & \to & K_2(\eta_X) &
\mapright{\oplus (\phantom{x},\phantom{x})_x} &
\bigoplus_{x\in X} \C^\ast \cr
& & & & \mapdown{c_2} & & \parallel \cr
0 & \to & H^1(X,\C^\ast) & \to & H^1(\eta_X,\C^\ast) &
\mapright{\oplus \delta_x}&  \bigoplus_{x\in X}\C\ast\cr }
$$
commutes. Since the rows are exact, the Chern class
$c_2$ induces a map $K_2(X) \to H^1(X,\C^\ast)$ which must be the
Chern class by naturality.

This construction extends easily to give a description of the
regulator
$$
c_2 : K_2(X) \to \Hdel^2(X,\Z(2))
$$
where $X$ is a smooth variety over $\C$. The construction proceeds in
the same way, except that the line bundles are no longer flat.
For a mixed Hodge structure $H$ denote
$\Hom_{\rm Hodge}(\Z,H)$, the set of ``Hodge classes'' in $H$ of type
$(0,0)$, by $\Gamma H$.

\begin{proposition}
If $X$ is smooth over $\C$, then there is a natural isomorphism
between $\Hdel^2(X,\Z(2))$ and the group which consists of the pairs
$(L,\nabla)$, where $L$ is a holomorphic line
bundles over $X$ and $\nabla$ is a holomorphic connection whose
curvature
times $2\pi i$ lies in $\Gamma H^2(X,\Z(2))$. \endproof
\end{proposition}

\section{The polylogarithm variation of mixed Hodge structure}
\label{polylog_varn}

Let $X$ be a smooth complex algebraic curve and $\X$ a smooth
compactification of it.  Let $D = \X -X$. Recall from
\cite{steenbrink-zucker}
that a {\it variation of mixed Hodge structure} over $X$ consists of
\begin{enumerate}
\item a $\Q$ local system $\V \to X$ which has a filtration by local
systems
$$
\cdots \subseteq \W_{l-1} \subseteq \W_l \subseteq \W_{l+1}
\subseteq \cdots
$$
which exhausts $\V$ and whose intersection is 0. We shall denote the
fiber
of $\V$ over $x\in X$ by $V_x$ and the fiber of $\W_l$ by $W_l V_x$.
We will
also assume that each local monodromy operator $T_P : V_P \to V_P$,
about
each $P \in D$, is unipotent.
\item a Hodge filtration
$$
\cdots \supseteq \cF^{p-1}\supseteq \cF^p \supseteq \cF^{p+1}
\supseteq \cdots
$$
of the corresponding holomorphic vector bundle $\cV := \V \otimes_\Q
\O_X$
by holomorphic sub-bundles. These are required to satisfy Griffiths'
transversality: If
$$
\nabla : \cV \to \cV \otimes_{\O_X} \Omega^1_X
$$
is the natural flat connection, then
$$
\nabla (\cF^p) \subseteq \cF^{p-1}\otimes_{\O_X} \Omega^1_X.
$$
Denote the fiber of $\cF^p$ over $x\in X$ by $F^p V_x$.
\item For each $x\in X$, the filtrations $W_\bullet V_x$ and
$F^\bullet V_x$ define a mixed Hodge structure on $V_x$.
\item Denote Deligne's canonical extension of $\cV$ to $\X$
by $\Vbar \to \X$ \cite{deligne_diffeq}. The Hodge bundles $\cF^p$
are
required to extend to holomorphic
sub-bundles $\Fbar^p$ of $\Vbar$. (Note that the weight bundles
$\cW_l := \W_l \otimes_\Q \O_X$ automatically extend to sub-bundles
$\Wbar$ of $\Vbar$ as they are flat.)
\item about each point $P\in D$, there is a relative weight
filtration
\cite{steenbrink-zucker}. This is an important condition which is
rather
technical in general. However, in the case where the global monodromy
representation
$$
\rho_x : \pi_1(X,x) \to GL(V_x)
$$
is unipotent, the condition reduces to the much simpler condition
$$
N_P(W_l V_x) \subseteq W_{l-2} V_x
$$
for each $P\in D$, where $N_P$ is the local monodromy logarithm
$$
N_P = {1 \over 2 \pi i} \log T_P.
$$
(See \cite{hain-zucker_1}.)
\end{enumerate}

\begin{theorem}\label{polylog-vmhs}
The $n$th polylogarithm local system underlies a good variation of
mixed
Hodge structure whose weight graded quotients are canonically
isomorphic
to $\Q,\Q(1), \ldots, \Q(n)$.
\end{theorem}

\proof
Let $\V \to \C - \01$ be the $n$th polylogarithm local
system, and $\cV$ the corresponding holomorphic vector bundle. By
(\ref{poly-loc-sys}), the canonical extension of this to $\P^1$ is
the
trivial bundle
$$
\P^1 \times \C^{n+1} \to \P^1.
$$
Denote the standard basis of $\C^{n+1}$ by
$e_0,e_1,\ldots, e_n$. The fiber $V_x$ is the $\Q$ linear span of
$\bL_0(x), \bL_1(x),\ldots, \bL_n(x)$, the rows of $\Lambda(x)$.
Define
\begin{equation}\label{filt_1}
W_{-2l + 1}\, \C^{n+1} = W_{-2l}\, \C^{n+1} =
\span\left\{e_l,\ldots,e_n\right\}
\end{equation}
and
\begin{equation}\label{filt_2}
F^{-p}\C^{n+1} = \span\left\{e_0,\ldots,e_p\right\}.
\end{equation}
Define
$$
\Fbar^p = \P^1 \times F^p\C^{n+1} \subseteq \Vbar
\hbox{ and } \Wbar_l = \P^1 \times W_l\C^{n+1} \subseteq \Vbar.
$$
Observe that the weight filtration comes from a filtration defined
on $\V$:
$$
W_{-2l + 1}V_x = W_{-2l}V_x =
\span\left\{\bL_l(x), \ldots, \bL_n(x)\right\}\quad \square
$$

Suppose that $\V \to X$ is a good variation of mixed Hodge structure
with unipotent monodromy about each point of $D = \X - X$. Let
$P \in D$. For each non-zero tangent vector $\v$ of $X$ at $p$, there
is a canonical mixed Hodge structure on $V_\C$, the fiber of $\Vbar$
over $P$. This is called the {\it limit mixed Hodge structure
associated
to} $\v$.  The Hodge and weight filtrations on $\V_\C$ are defined by
letting $F^p V_\C$ and $W_l V_\C$ be the fibers of $\Fbar^p$ and
$\Wbar_l$ over $P$, respectively. To construct the limit mixed Hodge
structure, we have to construct a rational form $V_\Q$ of $V_\C$
and show that the weight filtration defined above is the
complexification of a filtration of $V_\Q$.

To construct $V_\Q$, choose an imbedded closed disk $\Deltabar$ in
$\X$ centered at $P$. Let $t$ be a holomorphic parameter in
$\Deltabar$ such that $t(P) = 0$ and $|t|=1$ is $\partial \Deltabar$.
By choosing the disk to be small enough, we may suppose that
$\Deltabar -\{0\} \subseteq X$.

We first consider the case where $\v = \partial/\partial t$. Let
$x \in \Deltabar$ be the point corresponding to $t=1$. Choose a
$\Q$ basis $v_1,\ldots,v_m$ of $V_x$, the fiber of $\V$ over $x$.
Let $v_1(t),\ldots,v_m(t)$ be flat (possibly multivalued) sections
of $\V$ over $\Deltabar^\ast$ which satisfy $v_j(1) = v_j$ for each
$j$. Let $T : V_x \to V_x$ be the local monodromy operator, and
$N= \log T/2 \pi i$ be the local monodromy logarithm. For each $j$,
define
$$
s_j(t) =  t^{-N} v_j(t).
$$
Then each $s_j(t)$ is a single valued section of $\Vbar$ over
$\Deltabar^\ast$. In fact, by the construction of the canonical
extension $\Vbar \to \X$, the $s_j$ comprise a local framing of
$\Vbar$
over $\Deltabar$. In particular, $s_1(0),\ldots,s_m(0)$ is a $\C$
basis of $V_\C$. Define the rational form $V_\Q$ of $V_\C$ which
corresponds to $\partial /\partial t$ to be the $\Q$-linear span of
$s_1(0),\ldots,s_m(0)$.  By choosing the basis $v_1,\ldots,v_m$ of
$V_x$ to be adapted to its weight filtration, one can easily show
that the weight filtration of $V_\C$ is the complexification of a
filtration of $V_\Q$. The $\Q$ structure on $V_\C$ which corresponds
to the tangent vector $\v = \lambda \partial/\partial t$ is defined
to be
$$
V_\Q(\v) = \lambda^N V_\Q.
$$
It is not difficult to show that this rational structure depends only
on the tangent vector $\v$, and not on the choice of the parameter
$t$.

If the weight graded quotients of $\V \to X$ are constant as
variations of
Hodge structure (e.g. the polylogarithm variations), it is not
difficult
to show that $((V_\Q(\v),W_\bullet),(V_\C, F^\bullet))$ is a mixed
Hodge structure, and that $N: V \to V(-1)$ is a morphism of mixed
Hodge structures.

The following result is due to Deligne \cite{deligne_letter} in the
case
$n=2$. It is a straight forward computation using the monodromy
computation (\ref{monod}) and the procedure described above.

\begin{theorem}\label{limit}
Let $z$ be the natural coordinate function on $\C-\01$.
The limit mixed Hodge structure on the $n$th polylogarithm
variation at the tangent vector $\partial /\partial z$
at 0 has rational structure spanned by the vectors
$$
\pmatrix{s_0 \cr s_1 \cr \vdots \cr s_n} = \left(
\begin{array}{c| c c c c}
1 & 0 & 0 & \cdots & 0 \\ \hline
0 &  2\pi i &  0 & \cdots & 0  \\
\vdots & 0 & \ddots &  & \vdots \\
0 & \cdots & \cdots & 0  & (2\pi i)^n\\
\end{array}\right)
\pmatrix{e_0 \cr e_1 \cr \vdots \cr e_n}
$$
The limit mixed Hodge structure associated with the tangent vector
$-\partial/\partial z$ at 1 has rational structure spanned by the
vectors
$$
\pmatrix{s_0 \cr s_1 \cr \vdots \cr s_n} = \left(
\begin{array}{c| c c c c}
1 & 0 & \zeta(2) & \cdots & \zeta(n) \\ \hline
0 &  2\pi i &  0 & \cdots & 0  \\
\vdots & 0 & \ddots &  & \vdots \\
0 & \cdots & \cdots & 0  & (2\pi i)^n\\
\end{array}\right)
\pmatrix{e_0 \cr e_1 \cr \vdots \cr e_n}
$$
where $\zeta(s)$ denotes the Riemann zeta function.
In both cases, the Hodge and weight filtrations are defined as in
(\ref{filt_1}) and (\ref{filt_2}).
\end{theorem}

\section{Mixed Hodge structure on $\pi_1$}
\label{fund_gp}

In this section we give the construction of the mixed Hodge structure
on the fundamental group $\pi_1(X,x)$ in the special case where $X$
is a smooth variety over $\C$ whose $H^1(X)$ is pure of weight 2.
This purity condition may be restated in several equivalent ways.
For example, $H^1(X)$ is pure of weight 2 if and only if one and
(and hence all) smooth completions $\X$ of $X$ have first Betti
number 0. In particular, all Zariski open subsets of a Grassmannians
have this property.

To put a mixed Hodge structure on $\pi_1(X,x)$, it is necessary to
linearize it. The first step is to replace the fundamental group
by its group algebra $\Q\pi_1(X,x)$. This object is not well enough
behaved, and we need to linearize it further, which we do by
completion.
Let $J$ be the augmentation ideal. That is, $J$ is the kernel of the
augmentation
$$
\epsilon : \Q\pi_1(X,x) \to \Q
$$
which takes each $g\in \pi_1(X,x)$ to 1. The powers of $J$ define a
topology on $\Q\pi_1(X,x)$, and we consider
$$
\Q\pi_1(X,x)\comp = \lim_\rightarrow \Q\pi_1(X,x)/J^m,
$$
the $J$-adic completion.

The completed group ring is a complete Hopf algebra. That is, it is a
topological algebra and has a continuous algebra homomorphism
$$
\Delta : \Q\pi_1(X,x)\comp \to
\Q\pi_1(X,x)\comp \comptensor \Q\pi_1(X,x)\comp
$$
Here $\comptensor$ denotes completed tensor product.
This algebra homomorphism is induced by the usual diagonal
$$
\Q\pi_1(X,x) \to \Q\pi_1(X,x) \otimes \Q\pi_1(X,x)
$$
which takes each $g\in \Q\pi_1(X,x)$ to $g\otimes g$.

The Malcev Lie algebra $\g(X,x)$ associated to $\pi_1(X,x)$ is, by
definition, the set of primitive elements
$$
\left\{ X\in  \Q\pi_1(X,x)\comp : \Delta X =
1 \otimes X + X\otimes 1 \right\}
$$
of $\Q\pi_1(X,x)\comp$.
This is a complete topological Lie algebra. The bracket of the
elements
$A,B$ of $\g$ is their commutator $AB - BA$, which is again
primitive.
The topology is induced from that of $\Q\pi_1(X,x)\comp$. The
completed
group ring may be recovered from its primitive elements as its
completed
enveloping algebra $U\comp \g(X,x)$ (see \cite[Appendix A]{quillen}).

A mixed Hodge structure on $\Q\pi_1(X,x)\comp$ is, by definition,
a compatible sequence of mixed Hodge structures
$$
\cdots \to \Q\pi_1(X,x)/J^3 \to \Q\pi_1(X,x)/J^2 \to
\Q\pi_1(X,x)/J \to \Q \to 0
$$
on the truncations of the group ring. Note that each of these is a
finite dimensional vector space.

If the product and diagonal of $\Q\pi_1(X,x)\comp$ are morphisms
of mixed Hodge structure, then, as $\g(X,x)$ is the kernel of the
reduced diagonal
$$
\overline{\Delta} : \Q\pi_1(X,x)\comp \to
\left[\Q\pi_1(X,x)\comp/\Q\right] \comptensor
\left[\Q\pi_1(X,x)\comp/\Q \right] \approx J\comptensor J,
$$
the Malcev Lie algebra inherits a mixed Hodge structure compatible
with
its Lie algebra structure. This mixed Hodge structure determines the
one on $\Q\pi_1(X,x)\comp$ as $\g(X,x)$ generates $\Q\pi_1(X,x)\comp$
topologically.

\begin{theorem}{\rm \cite{morgan},\cite{hain} }
If $(X,x)$ is a complex algebraic variety, then $\Q\pi_1(X,x)\comp$
and
$\g(X,x)$ have canonical mixed Hodge structures which are compatible
with their algebraic structures. \end{theorem}

It is important to note that if $H^1(X,\Q)$ is non-trivial,
then this mixed Hodge structure depends non-trivially on the
basepoint
$x\in X$ (cf. \cite[\S\S6,7]{geom}).

We will sketch the construction of this mixed Hodge structure in the
case
when $H^1(X)$ is pure of weight 2. Write
$X = \X - D$, where $\X$ is smooth and complete, and $D$ is a divisor
with normal crossings in $\X$.  For convenience, we set
$$
\Omega^\bullet(X) = H^0(\X, \Omega_{\X}^\bullet(\log D)).
$$
By results of Deligne \cite[(3.2.14)]{deligne_II}, every element of
$\Omega^\bullet(X)$ is closed and the obvious map
$$
\Omega^\bullet(X) \to H^\bullet(X,\C)
$$
is injective. The tensor algebra
$$
T := \bigoplus_{n\in \N} \Omega^1(X)^{\otimes(-n)}
$$
on the dual of $\Omega^1(X)$
is a Hopf algebra; the diagonal is defined as the algebra
homomorphism
which takes each $X \in \Omega^1(X)^\ast$ to $1\otimes X + X\otimes
1$.
The ideal $(\im \delta)$ in $T$ generated by the image
of the dual of the cup product
$$
\delta : \Omega^2(X)^\ast \to
\Lambda^2 \Omega^1(X)^\ast \subseteq \Omega^1(X)^{\otimes(-2)}
$$
is a Hopf ideal. It follows that
$$
A :=  T/(\im \delta)
$$
is a Hopf algebra whose diagonal is induced from that of the tensor
algebra.  The set of primitive elements of $A$ is
$$
PA = \L(\Omega^1(X)^\ast)/(\im \delta),
$$
where $\L(V)$ denotes the free Lie algebra generated by the vector
space
$V$.

Denote the ideal generated by $\Omega^1(X)^\ast$ by $I$. The
powers of $I$ define a topology on $A$. Denote the $I$-adic
completion
of $A$ by $A\comp$. The set of primitive elements of $A\comp$ is the
$I$-adic completion of $PA$.

The following result is a special case of a theorem of K.-T. Chen
\cite[(3.5)]{chen}.

\begin{proposition}
For each $x\in X$ there is a canonical isomorphism
$$
\Theta_x : \C\pi_1(X,x)\comp \to A\comp
$$
of complete Hopf algebras.
\end{proposition}

\proof Let $\omega \in \Omega^1(X) \otimes \Omega^1(X)^\ast$ be the
element which corresponds to the identity
$\Omega^1(X) \to \Omega^1(X)$. This can be viewed as an element of
$\Omega^1(X)\otimes PA\comp$. This form is integrable. That is,
$$
d\omega + \omega \wedge \omega = 0.
$$
It follows that the value of the $A$-valued iterated integral
$$
1 + \int\omega + \int \omega\omega + \int\omega\omega\omega + \cdots
$$
on each path in $X$ depends only on its homotopy class relative to
its end points (cf. \cite[\S3]{geom}, for example). It follows from
this and
(\ref{props})(ii) that
this map induces a well defined homomorphism from $\pi_1(X,x)$ into
the group of units of $A\comp$. This extends to an algebra
homomorphism
$$
\C\pi_1(X,x) \to A\comp.
$$
Since the augmentation ideal of $\C\pi_1(X,x)$ is mapped into the
ideal
$I$ of $A\comp$, it follows that this homomorphism extends to a
continuous algebra homomorphism
$$
\Theta_x : \C\pi_1(X,x)\comp \to A\comp.
$$
The property  (\ref{props})(iv)  implies that $\Theta_x$ commutes
with
the diagonals; that is, $\Theta_x$ is a Hopf algebra homomorphism.

The graded module associated to the filtration of
$\C\pi_1(X,x)\comp$ by powers of $J$ is generated by $J/J^2$, and
this
is isomorphic to $H_1(X)$. Similarly, the graded module associated to
the filtration of $A\comp$ by the powers of $I$ is generated by
$I/I^2$,
which is also isomorphic to $H_1(X)$. The map $\Theta_x$ induces an
isomorphism $J/J^2 \approx I/I^2$. Since both algebras are complete,
this implies that $\Theta_x$ is surjective.  One can show, without
too
much difficulty, that the map $J^2/J^3 \to I^2/I^3$ induced by
$\Theta_x$
is also an isomorphism (cf. \cite[(6.1)]{geom}).  It is then
relatively
straightforward to show that
$\Theta_x$ must be injective. The idea is that $A\comp$ contains no
other
relations other than those that are consequences of the quadratic
ones,
while $\C\pi_1(X,x)\comp$ has at least these relations. Since
$\Theta_x$
is well defined, it must be an isomorphism. \endproof
\medskip

The next step in constructing the mixed Hodge structure on
$\Q\pi_1(X,x)\comp$ is to define the Hodge and weight filtrations.
We do this by defining them on $A\comp$ and transferring them to
$\Q\pi_1(X,x)\comp$ via the isomorphism $\Theta_x$.

The ring $A$ is graded as the ideal $(\im \delta)$ is graded. Write
$$
A = \bigoplus_{n\in \N} A_n
$$
where $A_n$ is the image of $\Omega^1(X)^{\otimes(-n)}$ in $A$.
The assumption that $H^1(X)$ be pure of weight 2 implies that it
is has Hodge type $(1,1)$. Consequently,
$\Omega^1(X)^\ast\approx H_1(X)$ has Hodge type $(-1,-1)$. It is
therefore
natural to define Hodge and weight filtrations on $A$ by
$$
F^{-p} A = \bigoplus_{n\le p} A_n
$$
and
$$
W_{-m} A = \bigoplus_{n\ge m/2}A_n
$$
Now define Hodge and weight filtrations on $\Q\pi_1(X,x)/J^l$ by
transferring the Hodge and weight filtrations from $A/I^l$ via
the isomorphism
$$
\Theta_x : \Q\pi_1(X,x)/J^l \to A/I^l.
$$

This data defines a mixed Hodge structure on
$\Q\pi_1(X,x)/J^l$. To see this, first observe that
$$
Gr^W_m\Q\pi_1(X,x)/J^l = \cases{
 J^r/J^{r+1} & when $m = -2r$ and $0\le r < l$;\cr
0 & otherwise.\cr}
$$
The Hodge filtration induced on $Gr^W_{-2p}$ satisfies
$$
F^{-p}\left[J^p/J^{p+1}\right] = J^p/J^{p+1},
\quad F^{-p+1}\left[J^p/J^{p+1}\right] = 0
$$
when $0\le p <l$. Because the weight filtration is defined over $\Q$,
it follows that these filtrations define a mixed Hodge structure on
$\Q\pi_1(X,x)/J^l$ whose weight graded quotients are all of even
weight,
and where the $2p$th graded quotient is of type $(p,p)$. Since the
multiplication and comultiplication of $A$ preserve the
filtrations and are defined over $\Q$, they are morphisms of mixed
Hodge structure.  It follows that $\g(X,x)$, endowed with the induced
filtrations, is a mixed Hodge structure.

We conclude this section by relating this mixed Hodge structure to
unipotent variations of mixed Hodge structure. A variation of
mixed Hodge structure over a smooth variety $X$ is good if its
restriction
to every curve satisfies the conditions in \S\ref{polylog_varn}. A
good
variation of mixed Hodge structure $\V \to X$ over a smooth variety
$X$
is {\it unipotent } if one (and hence all) monodromy representations
\begin{equation}\label{rep}
\rho_x : \pi_1(X,x) \to \Aut V_x
\end{equation}
are unipotent. This condition is equivalent to the condition that
each
of the variations of Hodge structure $Gr^W_m \V$ be constant.

The monodromy representation (\ref{rep}) induces a map
$$
\theta_x : \Q\pi_1(X,x) \to \End V_x.
$$
Since the representation is unipotent, there exists $l$ such that
$J^l$ is contained in $\ker \theta_x$. It follows that there is
an algebra homomorphism
\begin{equation}\label{ind_rep}
\theta_x : \Q\pi_1(X,x)/J^l \to \End V_x.
\end{equation}
Both sides of this last equation have natural mixed Hodge structures.

\begin{theorem}{\rm \cite{hain-zucker_1}}\label{morphism}
For each $ x\in X$, the representation (\ref{ind_rep}) is a morphism
of mixed Hodge structures.
\end{theorem}

Define the category of Hodge theoretic representations of
$\pi_1(X,x)$
to be the set of pairs $(V,\rho)$, where $V$ is a mixed Hodge
structure
and $\rho$ is a unipotent representation $\pi_1(X,x) \to \Aut V$
which induces a morphism of mixed Hodge structure
$$
\theta_x : \Q\pi_1(X,x)/J^l \to \End V
$$
for $l$ sufficiently large.
Theorem \ref{morphism} implies that taking the fiber at $x$ defines
a functor from the category of unipotent variations of mixed Hodge
structure over $X$ to the category of Hodge theoretic representations
of $\pi_1(X,x)$.

\begin{theorem}{\rm \cite{hain-zucker_1}}\label{equivalence}
This functor is an equivalence of categories.
\end{theorem}

The proofs of Theorems \ref{morphism} and \ref{equivalence} in the
case
when $H^1(X)$ is pure of weight 2  are considerably simpler than
in general case.  (See \cite{hain-zucker_2} for a proof in this
case.)

A good variation of mixed Hodge structure $\V$ over a smooth variety
$X$ is called a {\it Tate variation of mixed Hodge structure} if all
of its weight graded quotients are constant and of even weight, and
if each of the variations $Gr^W_{2p}\V$ is of type $(p,p)$. The
polylogarithm
variations are examples of Tate variations of mixed Hodge structure.
Now suppose that $\X$ is any smooth compactification of $X$ where
$D = \X - X$ is a divisor with normal crossings in $\X$. Let
$\Vbar \to \X$ be the canonical extension of $\V$ to $\X$.  The
following
result is a simple consequence of Theorem \ref{morphism} and
\cite[(6.4)]{riem_hilb}.

\begin{theorem}\label{trivial}
If $\V \to X$ is a Tate variation of mixed Hodge structure, then its
canonical extension $\Vbar \to \X$ is trivial as a holomorphic vector
bundle, so that there is a complex vector space $V$ and a bundle
isomorphism
$$
\matrix{
\Vbar & \to & V\times \X \cr
\downarrow & & \downarrow \cr
\X & = & \X \cr}
$$
Moreover, there are filtrations $F^\bullet$ and $W_\bullet$ of $V$
such that the extended Hodge and weight bundles $\Fbar^p$ and
$\Wbar_l$ correspond to $F^p \times \X$ and $W_l \times \X$,
respectively, under the bundle isomorphism. \endproof
\end{theorem}

One important example of a unipotent variation of mixed Hodge
structure over a smooth variety $X$ is the one whose fiber over
$x\in X$ is the truncated group ring $\Q\pi_1(X,x)/J^l$. This
is a good variation because the monodromy representation
$$
\Q\pi_1(X,x)/J^l \to \Aut \Q\pi_1(X,x)/J^l
$$
can be written in terms of left and right multiplication, and
is thus a morphism of mixed Hodge structure. Such variations form
an inverse system of variations, and we call the inverse limit
the {\it tautological variation} over $X$. In case when $H^1(X)$
has weight 2, this variation is a Tate variation of mixed Hodge
structure, and can be described explicitly. View $A\comp$ as a
subalgebra of $\End\, A\comp$ via the right regular representation.

\begin{proposition}\label{tautological}
If $H^1(X)$ is of weight 2, then the tautological variation over $X$
has canonical extension $A\comp \times \X \to \X$. The connection
form
of the canonical flat connection on this bundle is given by the
$PA\comp$ valued 1-form
$$
\omega \in \Omega^1(X) \otimes H_1(X) \subseteq
\Omega^1(X) \otimes PA\comp
\subseteq \Omega^1(X)\otimes \End\, A\comp
$$
which corresponds to the canonical isomorphism $\Omega^1(X)
\approx H^1(X)$. The extended Hodge and weight bundles
are $F^p A\comp \times \X$ and $W_l A\comp \times \X$.
\endproof
\end{proposition}

\section{Hodge theoretic interpretation of regulators}
\label{hodge-theory}

In this section we give a Hodge theoretic interpretation of the
regulators constructed in Sections \ref{reg-k3} and \ref{heisenberg}.
These interpretations are due to Deligne \cite{deligne_letter}.

Throughout this section, $\Lambda$ will denote $\Z,\Q$ or $\R$.
Suppose that $V=(V_\Lambda,(V_\LQ,W_\bullet),(V_\C,F^\bullet))$ is a
mixed Hodge structure where the underlying lattice $V_\Lambda$ is
torsion free. The ring
of endomorphisms $\End V$ has a mixed Hodge structure whose Hodge and
weight filtrations are defined by
$$
F^p \End_\C V =
\left\{ \phi \in \End_\C V : \phi(F^q V)\subseteq F^{p+q}V\right\}
$$
and
$$
W_l \End_\Q V =
\left\{ \phi \in \End_\LQ V : \phi(W_m V)\subseteq W_{m+l}V\right\}.
$$
Set $\g = W_{-1}\End V$. This is a nilpotent Lie algebra with a mixed
Hodge structure---the bracket being the commutator
$[\phi,\psi] = \phi\psi - \psi\phi$. The subspace $F^0\g$ is a Lie
sub-algebra. Denote the  simply connected Lie groups which correspond
to $\g_\C$ and $F^0\g$ by $G$ and $F^0 G$, respectively. These are
unipotent subgroups of $\Aut_\C V$. Set
$G_\Lambda = G \cap \End_\Lambda V$. We view $G_\Lambda$ as acting on
the right of $V_\Lambda$.

For each $g\in G$, the triple
$(V_\Lambda g,(V_\LQ g,W_\bullet g),(V_\C,F^\bullet))$
is a mixed Hodge structure whose weight graded quotients are
canonically
isomorphic to those of $V$. It is not difficult to show that every
mixed
Hodge structure with torsion free lattice and weight
graded quotients canonically isomorphic to those of $V$
can be constructed this way. More generally, we have the following
result
which is easily proved (cf. \cite{carlson}).

\begin{proposition}\label{moduli}
The set of $\Lambda$-mixed Hodge
structures whose weight graded quotients are canonically isomorphic
those of $V$ is naturally isomorphic to
$$
G_\Lambda\backslash G/F^0G.
$$
The double coset of $g\in G$ corresponds to the mixed Hodge structure
$$
V=(V_\Lambda g,(V_\LQ g,W_\bullet g),(V_\C,F^\bullet)).\quad \square
$$
\end{proposition}

This identification can be used to compute extension groups of mixed
Hodge
structures. Suppose that $A$ and $B$ are $\Lambda$-Hodge structures
whose underlying $\Lambda$ module is torsion free. Suppose that the
weight of $A$ is greater than that of $B$. If we take $V = A\oplus B$
then, for $R = \Lambda,\C$,
$$
G_R = \Hom_R(A,B)
$$
In this case the moduli space of mixed Hodge structures with weight
graded
quotients canonically isomorphic to $A$ and $B$ is the group
$\Ext_\cH^1(A,B)$
of extensions of $A$ by $B$ in the category $\cH$ of mixed Hodge
structures.
Applying Proposition (\ref{moduli}) to the split mixed Hodge
structure
$A\oplus B$, we obtain the well known formula for $\Ext_\cH^1$
(cf. \cite{carlson}).

\begin{proposition}\label{extensions}
With $A$ and $B$ as above, there is a canonical isomorphism
$$
\Ext_\cH^1(A,B) \approx
{\Hom_\C(A,B) \over \Hom_\Lambda(A,B) + F^0\Hom_\C(A,B)}.
$$
\end{proposition}

An important special case is where $A=\Z$, $B = \Z(n)$ and $n\ge 1$.
(Recall that $\Z(n)$ is the Hodge structure of type $(-n,-n)$ whose
lattice
is the subgroup $(2\pi i)^n\Z$ of $\C$.) In this case we have
$$
\Ext_\cH^1(\Z,\Z(n)) \approx \C/\Z(n).
$$
Following through the construction, we see that
the mixed Hodge structure which corresponds to $\lambda \in \C/\Z(n)$
can be described as follows. Denote the standard basis of $\C^2$ by
$e_0,e_n$. These have type $(0,0),(-n,-n)$, respectively. The
Hodge and weight filtrations on $\C^2$ are defined by
$$
W_l \C^2 = \span \{e_j: -j  \le l \}
$$
and
$$
F^p \C^2 = \span \{ e_j : -j \ge p\}.
$$
 The mixed Hodge structure which corresponds to $\lambda$
has integral basis the two vectors
$$
\pmatrix{1 & \lambda \cr 0 & (2\pi i)^n \cr}
\pmatrix{e_0 \cr e_n \cr }.
$$
In particular, the extension of $\Z$ by $\Z(1)$ which corresponds to
$x \in \C^\ast \approx \C/\Z(1)$ has integral basis spanned by the
vectors
$$
\pmatrix{1 & \log x \cr 0 & 2\pi i \cr}
\pmatrix{e_0 \cr e_1 \cr }.
$$

A unipotent variation of mixed Hodge structure over a smooth variety
$X$
whose weight graded quotients are canonically isomorphic to $\Z$ and
$\Z(m)$, $(m\ge 1)$ will determine a {\it classifying map}
$X \to \Ext^1_\cH(\Z,\Z(m))$.

\begin{proposition}\label{ext_varns}
When $m >1$ the classifying map is constant. When $m = 1$, a
map
$$
X \to \Ext^1_\cH(\Z,\Z(1)) \approx \C^\ast
$$
is the classifying map
of a good variation of mixed Hodge structure if and only if it is an
algebraic function on $X$.
\end{proposition}

\proof In both cases, the canonical extension of the variation to
a good compactification $\X$ of $X$ is a trivial bundle
(\ref{trivial}), as are the extended Hodge and weight bundles. In the
first case, Griffiths' transversality forces the integral lattice to
be constant. In the second, the regularity of the connection of the
canonical extension corresponds to the classifying map $X \to
\C^\ast$
of the variation having poles at infinity. \endproof

Since there is a canonical isomorphism
$$
\Ext_\cH^1(\Lambda,\Lambda(m)) \approx \C/\Lambda(m),
$$
the regulator $K_m(\C) \to \C/\Lambda(m)$
can then be interpreted as a map
$$
K_m(\C) \to \Ext_\cH^1(\Lambda,\Lambda(m)).
$$
A motivic description of this regulator in the case when $m=3$
is given in \cite{BGSV}.

The regulator
$$
c_2 : K_2(X) \to \Hdel^2(X,\Z(2))
$$
also admits a Hodge theoretic interpretation.  This time we take our
reference Hodge structure $V$ to be the direct sum of $\Z(0)$,
$\Z(1)$
and $\Z(2)$. The moduli space of mixed Hodge structures whose weight
graded quotients are canonically isomorphic to these Hodge structures
is
$$
H_\Z\backslash H_\C.
$$
where $H$ denotes the Heisenberg group defined in \S\ref{heisenberg}.
The bundle projection
$$
H_\Z \backslash H_\C \to \C^\ast \times \C^\ast
$$
may be interpreted as the map which takes a mixed Hodge structure
$V \in H_\Z\backslash H_\C$ to
$$
(V/\Z(2),W_2V)\in \Ext^1_\cH(\Z,\Z(1)) \times \Ext^1_\cH(\Z(1),\Z(2))
\approx \C^\ast \times \C^\ast.
$$

We next consider the problem of determining which maps
$X \to H_\Z\backslash H_\C$ classify variations of mixed Hodge
structure.

\begin{proposition}\label{classifying}
A function $f : X \to H_\Z\backslash H_\C$ is the classifying map
of a variation of mixed Hodge structure over $X$ with weight graded
quotients canonically isomorphic to $\Z,\Z(1)$ and $\Z(2)$ if and
only if\begin{enumerate}
\item f is holomorphic;
\item the composite
$$
X \mapright{f} H_\Z\backslash H_\C \to \C^\ast \times \C^\ast
$$
of $f$ with the canonical projection is algebraic;
\item the map $f : X \to H_\Z\backslash H_\C$ is a flat section
of the bundle $H_\Z\backslash H_\C^\ast \to \C^\ast \times \C^\ast$.
\end{enumerate}
\end{proposition}

\proof The first statement corresponds to the fact that the
connection
on the bundle $\cV = \V\otimes_\Z \O_X$ is holomorphic. The second
follows from (\ref{ext_varns}) and the fact that if $\V$ is a
variation, then so are $\V/\Z(2)$ and $W_2\V$. The last condition
corresponds to Griffiths' transversality. One needs to use the fact
that the canonical extension of $\V$ to a good compactification $\X$
of $X$ is trivial, and that the extended Hodge and weight bundles are
also trivial (\ref{trivial}). \endproof

This result allows us to give an interpretation of the regulator
$$
c_2 : K_2(X) \to \Hdel^2(X,\Z(2))
$$
constructed in Section \ref{heisenberg}: If $f,g$ are invertible
functions on $X$, then $c_2(\{f,g\})$ is the obstruction to finding
a good variation of mixed Hodge structure $\V$ over $X$ with weight
graded quotients $\Z,\Z(1),\Z(2)$ and whose subquotients $\V/\Z(2)$
and $W_{-2}\V$ are classified by
$$
f : X \to  \C^\ast \approx \Ext^1_\cH(\Z,\Z(1)) \hbox{ and }
g : X \to \C^\ast \approx \Ext^1_\cH(\Z(1),\Z(2)).
$$

More on extensions of variations of mixed Hodge structure can be
found
in \cite{carlson-hain} and \cite{alg_cycles}.

\section{The polylogarithm quotient of
$\pi_1(\protect\P^1-\{0,1,\infty\}$}
\label{polog_quot}

The polylogarithm quotient of the the fundamental group of
$\P^1-\{0,1,\infty\}$ is the image of the monodromy representation
$$
\pi_1(\P^1-\{0,1,\infty\},x) \to \Aut P_x
$$
where $P \to \P^1-\{0,1,\infty\}$ is the polylogarithm variation of
mixed Hodge structure.

Since the monodromy representation of the $n$th polylogarithm local
system is unipotent, it induces a representation
$$
\Q\pi_1(\C - \01,x)/J^{n+1} \to gl_{n+1}(\C).
$$
Denote the image of the composite
$$
\g(X,x) \to \Q\pi_1(\C - \01,x)/J^{n+1} \to gl_{n+1}(\C)
$$
by $\p_n(x)$.

\begin{proposition}
For each $x\in \C-\01$, $\p_n(x)$ has a natural mixed Hodge structure
compatible with its Lie algebra structure. The local system of the
$\p_n(x)$ forms a good unipotent variation of mixed Hodge structure
over
$\C-\01$.  Finally, these local systems form an inverse system of
variations of mixed Hodge structure.
\end{proposition}

\proof The first assertion is an immediate consequence of
(\ref{polylog-vmhs}) and (\ref{morphism}). The second is a
consequence
of (\ref{tautological}). The last assertion is clear. \endproof

By the construction given in \S\ref{fund_gp}, the Malcev Lie algebra
of $\pi_1(\C-\01,x)$ is the completion of the free Lie algebra
generated
by $H_1(\C-\01,\C)\approx \Omega^1(\C-\01)^\ast$. Let $X_0,X_1$ be
the
basis of $H_1(\C-\01,\Z)$ consisting of the homology classes of the
loops $\sigma_0$, $\sigma_1$ defined in \S\ref{monodromy}. This is
dual to the basis $\omega_0/2 \pi i$,
$- \omega_1/2\pi i$ of $\Omega^1(\C-\01)$, where $\omega_0,\omega_1$
are the forms defined in
\S\ref{monodromy}. The completed group ring $\C\pi_1(\C-\01,x)\comp$
is isomorphic to the completion $\C\langle\langle X_0,X_1
\rangle\rangle$
of the free associative algebra generated by $X_0,X_1$.  The set of
primitive elements of $\C\langle\langle X_0,X_1 \rangle\rangle$ is
$\f = \L(X_0,X_1)\comp$, the completion of the free Lie algebra
generated by $X_0$ and $X_1$. The isomorphism
$$
\C\pi_1(\C-\01,x)\comp \to \C\langle\langle X_0,X_1 \rangle\rangle
$$
is induced by the map $\pi_1(\C-\01,x)$ which takes $\gamma$ to
$$
1 + \int_\gamma \omega + \int_\gamma \omega\omega +
\int_\gamma \omega\omega\omega + \cdots
$$
where $\omega$ is the $\f$-valued 1-form
\begin{equation}\label{con_form}
\omega = \omega_0 X_0 - \omega_1 X_1.
\end{equation}

\begin{proposition}
The monodromy representation
$$
\pi_1(\C-\01,x)\to gl_{n+1}(\C)
$$
of the polylogarithm local system is induced by the homomorphism
$\f \to gl_{n+1}$ defined by
$$
X_0 \mapsto
\pmatrix{
0 & 0 & 0 & \cdots & 0 \cr
  & 0 & 1 & \ddots & \vdots \cr
\vdots & & \ddots & \ddots & 0 \cr
  & & & \ddots & 1 \cr
0 & & \cdots & & 0 \cr }
\quad X_1 \mapsto
\left(\begin{array}{c| r c c c}
0 & -1 & 0 & \cdots & 0 \\ \hline
0 & 0  &   0  & \cdots  & 0  \\
0 & 0 & 0 &  & 0 \\
\vdots & \vdots &  & \ddots & \\
0 & 0 & 0 & \cdots  & 0\\
\end{array}\right)
$$
\end{proposition}

\proof This follows as the connection matrix of the polylogarithm
local
system is
$$
\pmatrix{
0 & \w_1 & 0 & \cdots & 0 \cr
  & \ddots & \w_0 & \ddots & \vdots \cr
\vdots &  & \ddots & \ddots & 0 \cr
  &	&	& \ddots & \w_0 \cr
0 &	& \cdots &	& 0 \cr }
$$
which is simply the $\f$-valued form $\omega$ (\ref{con_form})
composed the the homomorphism
$\f \to gl_{n+1}(\C)$ defined in the statement of the proposition.
\endproof

\begin{corollary}
The complex form of the polylogarithm quotient has presentation
$$
\p = \L(X_0,X_1)\comp /(\ad(X_1)[\L,\L])
$$
and a topological basis $\left\{\ad (X_0)^n(X_1): n\in \N\right\}$.
\endproof
\end{corollary}

\section{Motivic Description of the Polylogarithm Variation}

In this section we give a motivic description of the polylogarithm
variations. This description goes back to Deligne.

First suppose that $X$ is a topological space. Denote the space of
paths
$\gamma : [0,1] \to X$ by $PX$. There is a canonical projection $PX
\to
X \times X$ which takes  a path $\gamma$ to its endpoints
$(\gamma(0),\gamma(1))$. Denote the fiber of this map over $(x,y)$ by
$P_{x,y}X$.

The group $H_0(P_{x,y}X)$ is the free abelian group on the set of
homotopy
classes of paths in $X$ from $x$ to $y$. These form a local system
\begin{equation}\label{canon_sys}
\left\{H_0(P_{x,y}X)\right\}_{(x,y)} \to X \times X.
\end{equation}
When $x=y$ there is a canonical isomorphism
$H_0(P_{x,y}X;\Q)\approx \Q\pi_1(X,x)$.  This has a canonical
filtration
given by the powers of the augmentation ideal $J$. This filtration
extends
to a flat filtration of the local system (\ref{canon_sys}). Denote
the completion of $H_0(P_{x,y}X;\Q)$ in the corresponding topology
by $H_0(P_{x,y}X;\Q)\comp$.

\begin{theorem}\label{canon_varn}
{\rm \cite{hain-zucker_1}}
If $X$ is a smooth algebraic variety, the local system
$$
\left\{H_0(P_{x,y}X;\Q)\comp\right\}_{(x,y)} \to X \times X
$$
is a good variation of mixed Hodge structure whose fiber over $(x,x)$
is
the canonical mixed Hodge structure on $\Q\pi_1(X,x)\comp$.
\end{theorem}

We call this the {\it canonical variation of mixed Hodge structure}
associated to $X$. Although this construction appears to be outside
the domain of algebraic geometry, it can be made motivic. There are
several equivalent ways of doing this. One can be found
in \cite{deligne:line}; the other is a  construction in topology of a
cosimplicial model of $PX$ and the fibration $PX \to X \times X$.
It is called the {\it cobar
construction} and dates back to the paper \cite{adams} of F.~Adams.
It makes sense for varieties over any base, as was noted by
Wojtkowiak
\cite{wojtkowiak}. Briefly, it is the cosimplicial space
$$
P^\bullet = X^{\Delta[1]_\bullet}
$$
where $\Delta[1]_\bullet$ is the standard simplicial model of the
unit interval. The projection $PX \to X\times X$ corresponds to the
map $X^{\Delta[1]_\bullet} \to X^{\partial \Delta[1]_\bullet}= X^2$
induced by the inclusion of the boundary of the interval.

Since the set of $m$-simplices of
$\Delta[1]_\bullet$ is the set of order preserving maps
$\{0,1,\ldots,m\} \to \{0,1\}$, and since there are precisely
$m+2$ of these,
$$
P^m = X \times X^m \times X.
$$
The coface maps $P^m \to P^{m+1}$ are the various diagonals. Applying
the de~Rham functor to this cosimplicial space yields the bar
construction
on the de~Rham complex of $X$, which is essentially  Chen's
complex of iterated integrals on $PX$ (\cite{chen},
see also \cite[\S\S 1,2]{hain}).

Suppose that $X = \X-D$ is a smooth curve and that $\v$ is a tangent
vector at a point $P\in D$. Suppose that $x\in X$. Denote by
$P_{\v,x}X$ the set of paths $\gamma :[0,1] \to \X$ which have the
property that $\gamma'(0) = \v$, $\gamma(1) = x$ and
$\gamma(]0,1]) \subseteq X$. This space is easily seen to be homotopy
equivalent to
$P_{z,x}X$ where $z$ is a point in $X$ which is sufficiently close
to $P$ in the direction of $\v$. More generally, one can define
$P_{\v_1,\v_2}X$ where $\v_1$ and $\v_2$ are non-zero tangent vectors
to points of $D$. Deligne defines $\pi_1(X,\v)$ to be the set of path
components of $P_{\v,\v}X$. It is canonically isomorphic to
$\pi_1(X,z)$
when $z\in X$ is sufficiently close to $P$ in the direction of $\v$
(cf.\ \cite{deligne:line}).

It is useful to think of the limit mixed Hodge structure of the
local system
$$
\left\{H_0(P_{z,x}X;\Q)\comp\right\}_{z\in X} \to X
$$
associated to $\v$ as a mixed Hodge structure on
$H_0(P_{\v,x}X;\Q)\comp$. Denote the limit mixed Hodge structure on
the Malcev Lie algebra of the polylogarithm quotient of
$\pi_1(\C-\01,\v)$ associated to the tangent vector $\v$ by $\p(\v)$.

\begin{theorem}
Let $\v$ be the tangent vector $\partial /\partial z$ at $0\in \P^1$.
The polylogarithm local system is the quotient of the variation of
mixed Hodge structure
$$
\left\{H_0(P_{\v,z}\C-\01;\Q)\comp\right\}_{z\in \C-\01}
\to \C-\01
$$
whose fiber over $\v$ is the polylogarithm quotient $\p(\v)$
of $\pi_1(\C-\01,\v)$.
\end{theorem}

Since both variations are isomorphic as local systems, to prove the
theorem it suffices to show that the fibers of of the two variations
over
one particular
point (or tangent vector) are isomorphic as mixed Hodge structures.
It is not difficult to show that the fiber of the
quotient of the canonical variation over the tangent vector
$\partial /\partial z$ at 0 is isomorphic to that of the polylog
variation,
which was calculated in (\ref{limit}).

\bibliographystyle{plain}

\begin{thebibliography}{99}

\bibitem{adams}
F.~Adams: On the cobar construction, {\it Colloque de Topologie
Alg\'ebrique (Louvain, 1956)}, George Thone, Liege, Masson,
Paris, 1957, 81--87.

\bibitem{beilinson}
A.~Beilinson:
{\it Higher regulators and values of $L$-functions},
 J.\ Soviet Math.\ 30 (1985), 2036--2070. Translated from:
Sovr.\ Probl.\ Mat.\ 24,  Mosc. VINITI (1984), 181--238.

\bibitem{BGSV}
A.~Beilinson, A.~Goncharov, V.~Schetmann, A.~Varchenko: { Aomoto
dilogarithms, mixed Hodge structures and motivic cohomology of pairs
of triangles in the plane}, in {\it The Grothendieck Festschrift :
a collection of articles written in honor of the 60th birthday of
Alexander Grothendieck}, P. Cartier et al., editors.
Birkhauser, Boston, 1990.


\bibitem{bloch_kyoto}
S.~Bloch: Applications of the dilogarithm function in algebraic
$K$-theory and algebraic geometry, {\it International Symposium on
Algebraic Geometry}, Kyoto, 1977, 103--114.

\bibitem{bloch-irvine}
S.~Bloch: {\it Higher regulators, Algebraic $K$-theory, and zeta
functions of elliptic curves}, unpublished manuscript, 1978.

\bibitem{borel}
A.~Borel: {\it Cohomologie de $SL_n$ et valeurs de fonctions de
zeta}, Ann.\ Scuola Normale Sup.\ 7 (1974), 613--636.

\bibitem{carlson}
J.~Carlson: {\it The geometry of the extension class of a mixed Hodge
structure}, Proc.\ Symp.\ Pure Math.\ 46, vol.\ 2 (1987), 199--222.

\bibitem{carlson-hain}
J.~Carlson, R.~Hain: {\it Extensions of variations of mixed Hodge
structure}, Asterisque 179--180 (1989), 39--65.

\bibitem{chen}
K.-T.~Chen: {\it Iterated integrals}, Bull.\ Amer.\ Math.\ Soc.\
83 (1977), 831--879.

\bibitem{deligne_diffeq}
P.~Deligne: {\it \'Equations Diff\'erentielles \`a  Points Singuliers
R\'eguliers}, LNM 163, Springer-Verlag, 1970.

\bibitem{deligne_II}
P.~Deligne: {\it Theorie d'Hodge, II}, Publ.\ Math.\ IHES 40 (1971),
5--58.

\bibitem{deligne_tame}
P.~Deligne:
Le symbole mod\'er\'e, unpublished notes, 1979.

\bibitem{deligne_letter}
P.~Deligne: Letter to Spencer Bloch, April 3, 1984.

\bibitem{deligne:line}
P.~Deligne: { Le groupe fondamental de la droite projective moins
trois points}, in {\it Galois groups over $\Q$ : proceedings of a
workshop held March 23-27, 1987}, editor K.~Ribet. Springer-Verlag,
New York, 1989.

\bibitem{dupont}
J.~Dupont: {\it The dilogarithm as a characteristic class for flat
bundles}, J.\ Pure and App.\ Alg.\ 44 (1987), 137--164.

\bibitem{D-H-Z}
J.~Dupont, R.~Hain, S.~Zucker: {\it Regulators and characteristic
classes of flat bundles}, preprint, January 1992.

\bibitem{dupont-sah}
J.~Dupont, C.-H.~Sah: {\it Scissors congruences, II}, J.\ Pure and
App.\ Alg.\ 24 (1982), 159--195.

\bibitem{goncharov}
A.~Goncharov: {\it Geometry of configurations, polylogarithms and
motivic
cohomology}, MPI preprint, 1991.

\bibitem{riem_hilb}
R.~Hain:
{\it On a generalization of Hilbert's 21st problem}, Ann.\ Sci.\
Ecole Norm.\ Sup.\ 19 (1986), 609--627.

\bibitem{geom}
R.~Hain: {\it The geometry of the mixed Hodge structure in the
fundamental
group}, Proc.\ Symp.\ Pure Math.\  46, vol.\ 2 (1987), 247--282.

\bibitem{hain}
R.~Hain: {\it The de~Rham homotopy theory of complex algebraic
varieties, I},
$K$-Theory 1 (1987), 271--324.

\bibitem{alg_cycles}
R.~Hain: {\it Algebraic cycles and extensions of variations of mixed
Hodge
structure}, Proc.\ Symp.\ Pure Math.\ 53 (1991), 175--221.

\bibitem{hain-macpherson}
R.~Hain, R.~MacPherson: {\it Higher logarithms}, Ill.\ J.\ Math.\ 34
(1990),
392--475.

preparation.

\bibitem{hain-zucker_1}
R.~Hain, S.~Zucker: {\it Unipotent variations of mixed Hodge
structure},
Invent.\ Math.\ 88 (1987), 83--124.

\bibitem{hain-zucker_2}
R.~Hain, S.~Zucker: {\it A Guide to unipotent variations of mixed
Hodge
structure}, Hodge Theory (Proceedings of the U.S. Spain Workshop,
Sant
Cugat, Spain, 1985), LNM 1246,  Springer-Verlag, 1987.

\bibitem{lewin}
L.~Lewin: {\it Polylogarithms and Associated Functions},
North-Holland,
New York, 1981.

\bibitem{lewin_2}
L.~Lewin: {\it Structural Properties of Polylogarithms}, Mathematical
Surveys and Monographs vol.\ 37, Amer.\ Math.\ Soc.,  1991.

\bibitem{morgan}
J.~Morgan: {\it The algebraic topology of smooth algebraic
varieties},
Publ.\ Math.\ IHES 48 (1978), 137--204, correction, Publ.\ Math.\
IHES 64 (1986), 185.

\bibitem{milnor_book}
J.~Milnor: {\it Introduction to Algebraic $K$-Theory}, Annals of
Math.\
Studies 72, Princeton University Press, 1971.

\bibitem{milnor}
J.~Milnor: {\it Hyperbolic geometry: the first 150 years}, Bull.\
Amer.\
Math.\ Soc.\ 6 (1982), 9--24.

\bibitem{quillen}
D.~Quillen: {\it Rational homotopy theory}, Ann.\ Math.\ 90 (1969),
205--295.

\bibitem{ramakrishnan_heisenberg}
D.~Ramakrishnan: {\it A regulator for curves via the Heisenberg
group},
Bull.\ Amer.\ Math.\ Soc.\ 5 (1981), 191--195.

\bibitem{ramakrishnan_monod}
D.~Ramakrishnan: {\it On the monodromy of higher logarithms}, Proc.\
Amer.\ Math.\ Soc.\ 85 (1982), 596--599.

\bibitem{ramakrishnan_bw}
D.~Ramakrishnan:  Analogs of the Bloch-Wigner function for higher
polylogarithms, in {\it Applications of Algebraic $K$-Theory to
Algebraic
Geometry and Algebraic Number Theory, part I}, Contemp.\ Math.\ 55
(1986),
371--376.

\bibitem{ramakrishnan_survey}
D.~Ramakrishnan: {\it Regulators, algebraic cycles, and values of
$L$-%
functions}, Contemp.\ Math.\ 83 (1989), 183--310.

\bibitem{steenbrink-zucker}
J.~Steenbrink, S.~Zucker: {\it Variations of mixed Hodge structure
I},
Invent.\ Math.\ 80 (1985), 489--542.

\bibitem{suslin}
A.~Suslin: Homology of $GL_n$, characteristic classes and Milnor
$K$-theory, in {\it Algebraic $K$-Theory, Number Theory, Geometry
and Analysis, Proceedings Bielefeld 1982}, LNM 1046, 357--375.

\bibitem{suslin_K3ind}
A.~Suslin: paper on $K_3^{ind}$.

\bibitem{wojtkowiak}
Z.~Wojtkowiak: {\it Cosimplicial objects in algebraic geometry},
preprint.

\bibitem{yang_thesis}
J.~Yang: {\it Algebraic $K$-groups of number fields and the
Hain-MacPherson
trilogarithm}, Ph.D. Thesis, University of Washington, 1991.

\bibitem{yang}
J.~Yang: {\it On the real cohomology of arithmetic groups and the
rank
conjecture for number fields}, Ann.\ Sci.\ Ecole
Norm.\ Sup.\ , to appear.

\bibitem{zagier}
D.~Zagier:  {\it The Bloch-Wigner-Ramakrishnan polylogarithm
function},
 Math.\ Ann.\  286, (1990), 613--624.

\end{thebibliography}

\end{document}